\documentclass[11pt,reqno]{amsart}
\usepackage[utf8]{inputenc}
\usepackage{amsmath, amsthm, amssymb}
\usepackage[margin=1in]{geometry}
\usepackage{hyperref}
\usepackage{mathtools} 
\usepackage{tikz} 
\usepackage{quiver} 
\usepackage{caption} 
\usepackage{subcaption}

\newtheorem{theorem}{Theorem}
\numberwithin{theorem}{subsection}
\newtheorem{lemma}[theorem]{Lemma}
\newtheorem{proposition}[theorem]{Proposition}

\newtheorem{corollary}[theorem]{Corollary}

\theoremstyle{definition}
\newtheorem{definition}[theorem]{Definition}
\newtheorem{example}[theorem]{Example}
\theoremstyle{remark}
\newtheorem*{remark}{Remark}

\newcommand{\CC}{\mathbb{C}}
\newcommand{\ZZ}{\mathbb{Z}}

\newcommand{\dist}{\normalfont \text{dist}}
\newcommand{\dom}{\normalfont \text{dom}}

\title{Iterated Lusztig--Vogan Bijection and Distinguished Weights}
\author{George Cao}
\address{Department of Mathematics, MIT, Cambridge, MA 02139}
\email{georgec8@mit.edu}
\date{}

\begin{document}

\begin{abstract}

The distinguished weights form a subset of the weight lattice and are closely tied to the notion of $p$-cells. These weights are defined via iterations of the Lusztig--Vogan bijection. We prove that all distinguished weights exhibit an anti-symmetry under the composition of reversal and negation. We show that the distribution of these weights follows a polynomial asymptotic, with a leading coefficient relating to the telephone numbers. As an explicit computation, we determine all the distinguished weights for $n \leq 4$.
\end{abstract}

\maketitle

\section{Introduction}

The representation theory of $GL_n(\CC)$ has been well-studied and is well-known. In particular, the irreducible finite-dimensional representations can be classified by their highest weight, which is always a dominant weight. On the other hand, the representation theory of $GL_n(\overline{\mathbb{F}}_p)$ is much more unknown, being now in the realm of modular representation theory. 

One tool to study this representation theory is the Lusztig--Vogan bijection, which was originally conjectured by both G.\,Lusztig~\cite{lusztig} and D.\,Vogan~\cite{vogan} independently. In generality, it states the following: let $G$ be a connected reductive algebraic group. Then this bijection is between $\Lambda^+$, the set of dominant weights of $G$, and $\Omega$, the set of pairs $(C, V)$ such that $C \subset G$ is a unipotent conjugacy class and $V$ is an irreducible representation of the centralizer $Z_{G}(g)$ for some $g \in C$. In the case of $G=GL_n(\CC)$, this is a bijection between $\Lambda_n^+$, the set of weakly decreasing (non-increasing) integer sequences of length $n$, and $\Omega_n$, the set of pairs $(\alpha, \nu)$ where $\alpha$ is a partition of $n$ and $\nu$ is an integer sequence that is dominant with respect to $\alpha$, i.e.~if $\alpha=(\alpha_1, \dots, \alpha_s)$ and $\nu = (\nu_1, \dots, \nu_t)$ where $\alpha_1 \geq \cdots \geq \alpha_s$, then $s=t$ and for all $i=1,\dots,s-1$ such that $\alpha_i = \alpha_{i+1}$ then $\nu_i \geq \nu_{i+1}$. 

In~\cite{bezrukavnikov}, R.\,Bezrukavnikov proved this bijection in full generality with novel and non-constructive methods. Prior to this, P.\,Achar gave an algorithm in~\cite{achar} for the case $G=GL_n$ to construct this bijection explicitly; the direction $\Lambda_n^+ \to \Omega_n$ is quite uninvolved, but the opposite direction is sophisticated (in fact, P.\,Achar proved this for $GL_n$ before the bijection was established to be true in general). In~\cite{rush}, D.\,Rush describes a simplified version of Achar's algorithm, which can be stated completely combinatorially. We will use Rush's algorithm since it is very concrete, and we are only concerned with $G = GL_n$.

The Lusztig--Vogan bijection therefore gives us a deeper understanding of the representation theory of reductive groups over characteristic $0$. To better understand the representation theory of $GL_n$ over a field of positive characteristic $p$, we consider another object of study: the set of distinguished weights, a subset of the dominant weights. Specifically, the distinguished weights are exactly the dominant weights for which, after a finite number of iterations of the algorithm for the Lusztig--Vogan bijection and division by $p$, only zeroes remain. The distinguished weights are related to $p$-cells \cite{jensen} and hence related to the modular representation theory of $GL_n$; in particular, each $p$-cell contains exactly one distinguished weight, so there is a natural bijection between the two sets. In this paper, we study these distinguished weights obtained through the iterative procedure.

\subsection{Main results} We state some of the main results of this paper. Let $\Lambda_{n,\text{dist}}^+$ be the set of distinguished weights of length $n$. Then we have the following anti-symmetry, a simple result with a lengthy proof.

\begin{theorem}[Corollary~\ref{cor:antisymmetric}]
If $\lambda \in \Lambda_{n,\dist}^+$ where $\lambda = (\lambda_1, \dots, \lambda_n)$, then $\lambda_i = -\lambda_{n+1-i}$.
\end{theorem}

Next, if we define $\Lambda_{n,k}^+$ to be the set of distinguished weights which map to all zeroes in at most $k$ iterations of the Lusztig--Vogan bijection algorithm and division by $p$. We find a recursive formula for $|\Lambda_{n,k}^+|$ in Theorem~\ref{thm:lambda_nk_recursion}, which leads to the following asymptotic formula:

\begin{theorem}[Theorem~\ref{thm:asymptotic}]
We have
$$|\Lambda_{n,k}^+| \sim \frac{a_{\left\lfloor \frac{n+1}2 \right\rfloor}}{\left\lfloor \frac n2 \right\rfloor !} k^{\left\lfloor \frac n2 \right\rfloor},$$
where $a_i$ is a sequence defined by $a_0 = a_1 = 1$ and $a_i = a_{i-1} + (i-1)a_{i-2}$.
\end{theorem}

This sequence of $a_i$ is commonly known as the telephone numbers and appears in more natural contexts, such as counting the number of involutions of $S_n$, the number of matchings on $n$ nodes, or the total number of standard Young tableaux on $n$ cells (OEIS~\cite[A000085]{oeis}).

Note that $\Lambda_{n,k}^+$ is the set of distinguished weights which consist of entries that are at most the order of magnitude of $O(p^{k-1})$ with respect to $p$, since each iteration of the algorithm approximately decreases the orders of magnitude by a factor of $p$. This gives a bound on the number of $p$-cells intersecting the cube $[-c,c]^n \subset \ZZ^n$ for some $c \in O_p(p^{k-1})$. We also include a discussion of the distribution of these distinguished weights, which have a regular pattern outside of the diagonals in $\ZZ^{n}$.

\subsection{Structure of the paper} The paper is structured as follows. In Section~\ref{sec:background}, we review the theory of the Lusztig--Vogan bijection, including the combinatorial description of the algorithm due to Rush~\cite{rush} (see Section~\ref{subsec:algorithm}). In Section~\ref{subsec:additional-background} we also describe a modification of the Lusztig--Vogan bijection, which we denote as $LV_p$, which is adapted to the modular setting; furthermore, we describe one of the main objects of study in this paper, the distinguished weights, in terms of $LV_p$. In Section~\ref{sec:examples}, we fully classify the distinguished weights for the cases $n=2,3,$ and $4$. In Section~\ref{sec:main-results}, we prove some general facts about the algorithm and distinguished weights, leading up to one of the main results, which gives a recursive formula and asymptotics for $|\Lambda_{n,k}^+|$, the number of distinguished weights of length $n$ which take at most $k$ iterations to become zero (Theorem~~\ref{thm:asymptotic}). Furthermore, we state and prove another main result that all distinguished weights have anti-symmetry (Corollary~\ref{cor:antisymmetric}). Finally, we have a short discussion about the distribution of these distinguished weights.

\section{Background}\label{sec:background}

\subsection{Notation}\label{subsec:notation}
We introduce notation that is used throughout the paper. Let $\alpha \vdash n$ denote that $\alpha$ is a partition of $n$. We may also explicitly write the parts of the partition as $\alpha = (\alpha_1, \alpha_2, \dots, \alpha_s)$ to denote that $\alpha_i \in \ZZ$, $\sum_{i=1}^s \alpha_i = n$, and $\alpha_1 \geq \alpha_2 \geq \cdots \geq \alpha_s \geq 1$. Then we can write $|\alpha|:=s$ as the number of parts. We may also write $\alpha = (a_1^{m_1},a_2^{m_2},\dots, a_t^{m_t})$ where $a_1 > a_2 > \cdots > a_t \geq 1$ and $\sum_{i=1}^t a_i m_i = n$, which emphasizes the multiplicity of each part.

We denote $\dom(\cdot)$ to be the function mapping a finite integer sequence (or set) to the ordered permutation of itself which is sorted in weakly decreasing (i.e.~non-increasing) order.

Throughout this paper we fix $p$ to be a sufficiently large prime. We expect that this bound need not be particularly large; for example, we expect that $p>2h$ for $h$ the Coxeter number should be sufficient, but we will not concern ourselves with precise lower bounds. Of course, the ``distinguished" weights and ideals depend on $p$, but we will suppress this from our notation.

\subsection{Algorithm for Lusztig--Vogan bijection}\label{subsec:algorithm}
We first introduce the integer sequences algorithm for computing the Lusztig--Vogan bijection, keeping the notation in~\cite{wang}. The algorithm for the Lusztig--Vogan bijection is the composition of three maps: $\phi$, $E^{-1}$, and $\kappa$. First, we describe the spaces that these maps are between.

\begin{definition}
We define $\Lambda_n^+$ as the set of all weakly decreasing integer sequences of length $n$, i.e.~$\lambda \in \Lambda_n^+$ if $\lambda = (\lambda_1, \dots, \lambda_n)$ where $\lambda_i \in \ZZ$ for all $1 \leq i \leq n$ and $\lambda_1 \geq \cdots \geq \lambda_n$.
\end{definition}

\begin{definition}\label{def:omega_n}
We define $\Omega_n$ as the set of pairs $(\alpha, \nu)$, where $\alpha = (\alpha_1, \dots, \alpha_l)$ is a partition of $n$ and $\nu = (\nu_1, \dots, \nu_l)$ is a dominant integer sequence with respect to $\alpha$, i.e.~if $\alpha_i = \alpha_{i+1}$, then $\nu_i \geq \nu_{i+1}$ for all $i=1, \dots, l-1$. 

For the algorithm, we can equivalently define $\Omega_n$ to be the set of tuples $(\mu_1, \dots, \mu_s)$ such that $\mu_i$ are weakly decreasing integer sequences (possibly empty), $\mu_s \neq \emptyset$, and $\sum_{i=1}^s i |\mu_i| = n$. We will be using this definition throughout the paper.

To see that these definitions are equivalent, observe that if we have a pair $(\alpha, \nu)$ such that $\alpha \vdash n$ and $\nu$ is a dominant integer sequence with respect to $\alpha$, then let $\mu_i = \dom(\{\nu_j : \alpha_j = i\})$ for\\ $1 \leq i \leq |\alpha|$. On the other hand, given a tuple $(\mu_1, \dots, \mu_s)$ satisfying the above conditions, then let $\alpha = (s^{|\mu_s|},(s-1)^{|\mu_{s-1}|},\cdots, 1^{|\mu_1|})$ and $\nu$ to be the concatenation of $\dom(\mu_s), \dom(\mu_{s-1}), \dots, \dom(\mu_1)$, in that order.
\end{definition}

\begin{definition}
We call $X$ a \textit{weighted diagram of size $n$} if $X = (X_1, \dots, X_l)$ is a tuple of weakly decreasing sequences of integers, where $n = \sum_{i=1}^l |X_i|$. Denote $X_{i,j}$ as the $j$th element of $X_i$, for $1 \leq j \leq |X_i|$. Call $X_i$ as the $i$th row of $X$, and denote the $j$th column of $X$ as $\{X_{i,j} : 1 \leq i \leq l, 1 \leq j \leq |X_i|\}$. Let $D_n$ be the set of all weighted diagrams of size $n$.
\end{definition}

We can begin defining the maps between these spaces that are part of the algorithm. In order to introduce the first map, we give the following definition, as appears in~\cite{wang}.

\begin{definition}
Let $\sigma \in \Lambda_n^+$. The \textit{maximal clumping of $\sigma$} is the unique partition of $\sigma$ into sequences $A_1, \dots, A_l$, where $l$ is minimal such that the concatenation of the $A_1, \dots, A_l$, in that order, is $\sigma$, and the underlying set of $A_i$ consists of consecutive numbers. Call each $A_i$ a \textit{maximal clump} or just \textit{clump}.
\end{definition}

\begin{definition}[Construction of $\phi$]
Define $\phi: \Lambda_n^+ \to D_n$ by the following iterative construction:

Let $\sigma \in \Lambda_n^+$. We build the image $X = \phi(\sigma)$ column-by-column. To start, let $\sigma_1 = \sigma$. To build the $r$th column, suppose we are given $\sigma_r$, a weakly decreasing integer sequence. Let $\sigma_r = A_1 \cup \cdots \cup A_l$ be its maximal clumping, and let $B_i$ be the sequence $A_i$ with repetitions removed. We construct a set $Z_r$ as follows:

For $1 \leq i \leq l$, let $B_i = (b_1, b_2, \dots)$. If both $r$ and $|B_i|$ are even, then append to $Z_r$ the elements $b_2, b_4, b_6, \dots$. Otherwise, append to $Z_r$ the elements $b_1, b_3, b_5, \dots$.

If $r=1$, then the first column of $X$ is $Z_r$, arranged in decreasing order from top to bottom and in rows $1$ to $|Z_r|$. Otherwise, the $r$th column of $X$ is still $Z_r$, arranged in decreasing order from top to bottom, with the additional condition that for every $z \in Z_r$, it must be placed in row $i$ such that $X_{i, r-1}$ is $z$ or $z + (-1)^r$. Then remove the elements of $Z_r$ from $\sigma_r$ to create $\sigma_{r+1}$ (for repetitions, only remove from $\sigma_r$ the multiplicity of the element in $Z_r$).

Repeat this procedure until all elements are removed. The validity of this algorithm is proved in~\cite{wang}.
\end{definition}

\begin{remark}
There is clearly an inverse map $\phi^{-1}$, which sorts the elements of $X \in D_n$ in weakly decreasing order.
\end{remark}

\begin{definition}[Construction of $E$]
Define $E: D_n \to D_n$ by the following element-wise operation:

Let $X \in D_n$ with $X = (X_1, \dots, X_l)$. We construct $E(X) = (X_1', \dots, X_l')$ satisfying the following conditions. First, we have $|X_i| = |X_i'|$. Second, for every $X_{i,j} \in X$, we have
$$X_{i,j}' = X_{i,j} + 2m_{i,j} - (c_j-1),$$
where $c_j = |\{r : 1 \leq r \leq l, 1 \leq j \leq |X_r|\}|$ is the total number of elements in column $j$ and $m_{i,j}$ is the number of $r$ such that $X_{r,j} \in X$ and $(X_{r,j}, -r)$ is lexicographically before $(X_{i,j}, -i)$.
\end{definition}

\begin{remark}[Construction of $E^{-1}$]
If $X$ is in the image of $\phi$, then each column of $X \in D_n$ is strictly decreasing with consecutive differences at least $2$ by construction. In this case, there is an inverse map $E^{-1}$, constructed as follows: 

Let $X \in D_n$ with $X = (X_1, \dots, X_l)$ satisfying the above properties, i.e.~$X_{r_1,j} - X_{r_2,j} \geq 2$ if $r_1 < r_2$. Then $E^{-1}(X) = (X_1', \dots, X_l')$, satisfying the following. First, $|X_i| = |X_i'|$. Second, for every $X_{i,j} \in X$, we have $X_{i,j}' = X_{i,j} + (c_j - 1) - 2m_{i,j}$, where $c_j$ and $m_{i,j}$ are defined the same as the above definition. Since the numbers in each column of $X$ are already sorted in decreasing order from top to bottom, this is the same as adding $(-(c_j-1), -(c_j-3), \dots, c_j-1)$ to the $j$th column, from top to bottom.
\end{remark}

\begin{definition}[Construction of $\kappa$]
Define $\kappa: D_n \to \Omega_n$ as follows:

Let $X \in D_n$ with $X = (X_1, \dots, X_l)$. Then $\kappa(X) = (\mu_1, \mu_2, \dots, \mu_s)$ where $s = \max_{1 \leq i \leq l}\{ |X_i| \}$ and $\mu_i = \dom(\{\sum_{j=1}^{|X_r|} X_{r,j}: |X_r| = i\})$, i.e.~the weakly decreasing sequence of row sums of rows with length $i$.
\end{definition}

\begin{remark}
The simplified construction for $\kappa^{-1}$ is the algorithm $\mathfrak{A}(\alpha, \nu)$ given in~\cite{rush}. It is omitted here due to its length.
\end{remark}

With these maps, we can define the map $LV : \Lambda_n^+ \to \Omega_n$ by $LV = \kappa \circ E^{-1} \circ \phi$. This is the algorithm for one direction of the Lusztig--Vogan bijection.

\begin{example}\label{ex:algorithm}
Let $n=8$ and let $\sigma = (46, 46, 45, 1, -1, -45, -46, -46) \in \Lambda_n^+$. Then
\begin{center}
    \begin{tikzpicture}[scale=0.8]
        \node at (-6.85,0) {$(46, 46, 45, 1, -1, -45, -46, -46) \xmapsto{\phi} $};
        \draw (-3,1) rectangle node{$46$} ++(1,1);
        \draw (-3,0) rectangle node{$1$} ++(1,1);
        \draw (-3,-1) rectangle node{$-1$} ++(1,1);
        \draw (-3,-2) rectangle node{$-45$} ++(1,1);
        \draw (-2,1) rectangle node{$45$} ++(1,1);
        \draw (-2,-2) rectangle node{$-46$} ++(1,1);
        \draw (-1,1) rectangle node{$46$} ++(1,1);
        \draw (-1,-2) rectangle node{$-46$} ++(1,1);
        \node at (0.75,0) {$\xmapsto{E^{-1}}$};
        \draw (1.5,1) rectangle node{$43$} ++(1,1);
        \draw (1.5,0) rectangle node{$0$} ++(1,1);
        \draw (1.5,-1) rectangle node{$0$} ++(1,1);
        \draw (1.5,-2) rectangle node{$-42$} ++(1,1);
        \draw (2.5,1) rectangle node{$44$} ++(1,1);
        \draw (2.5,-2) rectangle node{$-45$} ++(1,1);
        \draw (3.5,1) rectangle node{$45$} ++(1,1);
        \draw (3.5,-2) rectangle node{$-45$} ++(1,1);
        \node at (7.35,0) {$\xmapsto{\kappa} ((0,0), (), (132, -132))$};
    \end{tikzpicture}
\end{center}
Hence, $LV(\sigma) = ((0,0), (), (132, -132))$.
\end{example}

\subsection{Iterations of the algorithm}\label{subsec:additional-background}

Fix a prime $p > n$. Define $LV_p = \frac1p \circ LV$, where $\frac1p$ denotes dividing by $p$. In particular, if $LV(\sigma) = (\mu_1, \dots, \mu_s)$ for $\sigma \in \Lambda_n^+$, then $LV_p(\sigma) = (\frac{\mu_1}p, \dots, \frac{\mu_s}p)$ where $\frac{\mu_i}p$ denotes dividing each element in the sequence by $p$. Note that $LV_p$ may not always be defined, because division by $p$ may not always become integers. However, in this paper, we are only concerned with the subset of $\Lambda_n^+$ which this map is defined, i.e.~becomes integral after this division by $p$.

The output of $LV_p$, assuming it is integral, is an element of $\Omega_n$, which can be written as $(\mu_1, \dots, \mu_s)$, where each $\mu_i$ is a weakly decreasing sequence. This means that $\mu_i \in \Lambda_{n_i}^+$ for some $n_i$, so we can apply $LV_p$ again to each $\mu_i$, only if the result is still integral and if $n_i \geq 2$. Then we obtain a tuple of tuples of weakly decreasing sequences, so we can further iterate again. In this way, we can keep iterating $LV_p$ until every sequence is either of length $0$, length $1$, or another iteration would make the numbers non-integral.

\begin{definition}
Let $\Lambda_{n, \text{dist}}^+ \subset \Lambda_n^+$ be the set of \textit{distinguished weights} of $GL_n$; in regards with the above algorithm, these are exactly the elements in $\Lambda_n^+$ such that upon iterations of $LV_p$, every remaining number is a zero. Note that
$$LV_p: (0, \dots, 0) \mapsto ((), \dots, (), (0)),$$
so we can stop iterating on a sequence once it has become all zeroes.

Furthermore, define $\Lambda_{n,k}^+ \subset \Lambda_{n,\dist}^+$ as the weakly decreasing sequences of length $n$ which become all zeroes after at most $k$ iterations of the algorithm. It is clear that $\Lambda_{n,0}^+ = \{(0, \dots, 0)\}$ and $\Lambda_{n,0}^+ \subset \Lambda_{n,1}^+ \subset \Lambda_{n,2}^+ \subset \cdots$.
\end{definition}

The following example gives an example of a distinguished weight.

\begin{example}
Continuing Example~\ref{ex:algorithm} for $n=8$ and $p=11$, we find that
\begin{align*}
    \sigma = (46,46,45,1,-1,-45,-46,-46) &\xmapsto{LV_p} ((0,0),(),(12,-12))\\
    &\xmapsto{LV_p} ((0,0),(),(1,-1))\\
    &\xmapsto{LV_p} ((0,0),(),(0,0)).
\end{align*}
Since three iterations of $LV_p$ on $\sigma$ becomes all zeroes, then $\sigma \in \Lambda_{n,\text{dist}}^+$. In particular, we have that $\sigma \in \Lambda_{n,k}^+$ for $k \geq 3$.
\end{example}

The map $LV_p$ returns an element in $\Omega_n$, which can be written as $(\alpha^{(1)}, \nu)$, where $\alpha^{(1)}$ is a partition of $n$ and $\nu$ is an integer sequence dominant with respect to $\alpha^{(1)}$. Let $\alpha^{(1)} = (a_1^{m_1}, \dots, a_s^{m_s})$. Then when we iterate by applying $LV_p$ again, we obtain further partitions $\alpha_1^{(2)} \vdash m_1$, $\alpha_2^{(2)}\vdash m_2,$ and so on. This tuple of partitions $(\alpha_1^{(2)}, \dots, \alpha_s^{(2)})$ is a refinement of the partition $\alpha^{(1)}$. When we iterate further, we recursively obtain further refinements of these partitions. In this manner, each distinguished weight corresponds to some sequence $(\alpha^{(1)}, \alpha^{(2)}, \dots)$ where $\alpha^{(1)}$ is a partition of $n$ and each successive element is a refinement of the previous.

\section{Characterization of $\Lambda_{n,{\normalfont \text{dist}}}^+$ for small $n$}\label{sec:examples}

\begin{example}
Let $n=2$. Then $\Lambda_{n, \dist}^+ = \{(\frac{p^m-1}{p-1}, -\frac{p^m-1}{p-1}) : m \geq 0\}$. When iterating, we find that
$$\left(\frac{p^m-1}{p-1}, -\frac{p^m-1}{p-1}\right) \xmapsto{m \text{ iterations}} \underbrace{(\cdots(}_{m+1}0, 0) \cdots).$$
The algorithm works as follows:
{\normalsize
\begin{center}
    \begin{tikzpicture}[scale=0.82]
        \node at (-7.75,0) {$(1+p+\cdots+p^{m-1}, -(1+p+\cdots+p^{m-1})) \xmapsto{\phi} $};
        \draw (-3,0) rectangle node{$1+p+\cdots+p^{m-1}$} ++(4.5,1);
        \draw (-3,-1) rectangle node{$-(1+p+\cdots+p^{m-1})$} ++(4.5,1);
        \node at (2.5,0) {$\xmapsto{E^{-1}}$};
        \draw (3.25,0) rectangle node{$p+\cdots+p^{m-1}$} ++(4,1);
        \draw (3.25,-1) rectangle node{$-(p+\cdots+p^{m-1})$} ++(4,1);
    \end{tikzpicture}
\end{center}
$$\xmapsto{\kappa} ((p+\cdots+p^{m-1}, -(p+\cdots+p^{m-1}))) \xmapsto{\frac1p} ((1+p+\cdots+p^{m-2}, -(1+p+\cdots+p^{m-2}))),$$
}
which arrives to zeroes recursively. In particular, we have $|\Lambda_{2,k}^+| = k+1$.

We claim that these are the only distinguished weights for $n=2$. This is because a distinguished weight must end with all zeroes, and if $n=2$, it must end at $(\cdots (0,0) \cdots )$ or $(\cdots ((), (0)) \cdots )$, which correspond to the partitions $(1,1)$ and $(2)$ of $2$, respectively. The second case comes directly from the first after one iteration, since $LV_p : (0,0) \mapsto ((), (0))$. By the bijectivity of the algorithm, the characterization above consists of all possible distinguished weights for $n=2$.
\end{example}

\begin{example}
Let $n=3$. Then
$$\Lambda_{n,\dist}^+ = \left\{\left(2\cdot\frac{p^m-1}{p-1}, 0, -2\cdot \frac{p^m-1}{p-1}\right) : m \geq 0\right\} \cup \left\{\left(\frac{p^{m+1} + p^m - 2}{p-1}, 0, -\frac{p^{m+1} + p^m - 2}{p-1}\right): m \geq 0 \right\}.$$
When iterating, we find that
$$\left(2\cdot\frac{p^m-1}{p-1}, 0, -2\cdot \frac{p^m-1}{p-1}\right) \xmapsto{m \text{ iterations}} \underbrace{(\cdots(}_{m+1}0,0,0)\cdots),$$
and
$$\left(\frac{(p^{m+1}-1) + (p^m - 1)}{p-1}, 0, -\frac{(p^{m+1}-1) + (p^m - 1)}{p-1}\right)  \xmapsto{m+1 \text{ iterations}} \underbrace{\Big(\cdots\Big(}_{m+1}(0),(0)\Big)\cdots\Big).$$
We have that $|\Lambda_{3,k}^+| = 2k+1$.

We show that these are all of the distinguished weights for $n=3$ by a similar argument as the previous example. In particular, since the algorithm corresponds to refinements of partitions of $n$, then a distinguished weight must end in $( \cdots (0,0,0)\cdots)$, $(\cdots ((0),(0)) \cdots )$, or $(\cdots ((),(),(0))\cdots)$, which correspond to the partitions $(1,1,1), (1,2),$ and $(3)$ of $n=3$, respectively. The third case comes directly from the first one after one more iteration, and by the bijectivity of the algorithm, the characterization above for $\Lambda_{n, \text{dist}}^+$ covers all distinguished weights for $n=3$.
\end{example}

\begin{example}
Let $n=4$. We have the following families of distinguished weights:
$$\frac{p^m-1}{p-1} \cdot (3,1,-1,-3) \xmapsto{m \text{ iterations}}  \underbrace{(\cdots(}_{m+1}0,0,0, 0)\cdots),$$
and
$$\left(\frac{p^{m+1}+2p^m-3}{p-1}, \frac{p^m-1}{p-1}, -\frac{p^m-1}{p-1}, -\frac{p^{m+1}+2p^m-3}{p-1}\right) \xmapsto{m+1 \text{ iterations}} \underbrace{\Big(\cdots\Big(}_{m+1}(0), (), (0) \Big) \cdots \Big).$$
In addition, there are more complicated closed forms. One of them is:
\begin{align*}
    &\left(\frac{p^{m+k+1} + p^{k+1} + p^k - 3}{p-1}, \frac{p^k-1}{p-1}, -\frac{p^k-1}{p-1}, -\frac{p^{m+k+1} + p^{k+1} + p^k - 3}{p-1} \right) \\
    &\quad\xmapsto{k \text{ iterations}} \underbrace{\Big( \cdots \Big(}_{k+1} \frac{p^{m+1}+p-2}{p-1}, 0, 0, -\frac{p^{m+1}+p-2}{p-1}\Big) \cdots \Big) \\
    &\quad\xmapsto{1 \text{ iteration}} \underbrace{\Big( \cdots \Big(}_{k+1} \left(\frac{p^m-1}{p-1}, -\frac{p^m-1}{p-1}\right), (0)\Big) \cdots \Big) \\
    &\quad\xmapsto{m \text{ iterations}} \underbrace{\Big( \cdots \Big(}_{k+1} \underbrace{(\cdots(}_{m+1}0, 0) \cdots), (0)\Big) \cdots \Big).
\end{align*}
The next closed form is split into two cases by parity. If $m$ is even, then we have
\begin{align*}
    &\left(\frac{p^{m+k} + 2p^{k+1} + 3p^k - 6}{2(p-1)}, \frac{p^{m+k}+p^k-2}{2(p-1)}, -\frac{p^{m+k}+p^k-2}{2(p-1)},-\frac{p^{m+k} + 2p^{k+1} + 3p^k - 6}{2(p-1)} \right) \\
    &\quad\xmapsto{k \text{ iterations}} \underbrace{\Big( \cdots \Big(}_{k+1} \frac{p^m+2p-3}{2(p-1)}, \frac{p^m-1}{2(p-1)}, -\frac{p^m-1}{2(p-1)}, -\frac{p^m+2p-3}{2(p-1)}\Big) \cdots \Big) \\
    &\quad\xmapsto{1 \text{ iteration}} \underbrace{\Big( \cdots \Big(}_{k+1} (), \left(\frac{p^{m-1}-1}{p-1}, -\frac{p^{m-1}-1}{p-1}\right)\Big) \cdots \Big) \\
    &\quad\xmapsto{m-1 \text{ iterations}} \underbrace{\Big( \cdots \Big(}_{k+1} (), \underbrace{(\cdots(}_{m}0, 0) \cdots)\Big) \cdots \Big).
\end{align*}
If $m$ is odd, then we have
\begin{align*}
    &\left(\frac{p^{m+k} + p^{k+1} + 4p^k - 6}{2(p-1)}, \frac{p^{m+k}+p^{k+1}-2}{2(p-1)}, -\frac{p^{m+k}+p^{k+1}-2}{2(p-1)},-\frac{p^{m+k} + p^{k+1} + 4p^k - 6}{2(p-1)} \right) \\
    &\quad\xmapsto{k \text{ iterations}} \underbrace{\Big( \cdots \Big(}_{k+1} \frac{p^m+p-2}{2(p-1)}, \frac{p^m+p-2}{2(p-1)}, -\frac{p^m+p-2}{2(p-1)}, -\frac{p^m+p-2}{2(p-1)}\Big) \cdots \Big) \\
    &\quad\xmapsto{1 \text{ iteration}} \underbrace{\Big( \cdots \Big(}_{k+1} (), \left(\frac{p^{m-1}-1}{p-1}, -\frac{p^{m-1}-1}{p-1}\right)\Big) \cdots \Big) \\
    &\quad\xmapsto{m-1 \text{ iterations}} \underbrace{\Big( \cdots \Big(}_{k+1} (), \underbrace{(\cdots(}_{m}0, 0) \cdots)\Big) \cdots \Big).
\end{align*}
Carefully counting, we find that $|\Lambda_{4,k}^+| = k^2 + 3k + 1$.

We can show that these are all of the distinguished weights for $n=4$ by a similar argument as the previous two examples.
\end{example}

\section{Main results on distinguished weights}\label{sec:main-results}

\subsection{Basic results}\label{subsec:basic-results}

When computing the distinguished weights, we could start at the end sequence of zeroes in many nested parentheses and use $LV_p^{-1}$ to obtain the distinguished weight which has that particular end sequence. Since we know this algorithm is a bijection, then each end sequence unique corresponds to a distinguished weight. The following proposition gives an exact formula for iterating $LV_p^{-1}$ on a single weakly decreasing sequence which is nested in multiple parentheses.

\begin{proposition}\label{prop:one_sequence_iterate}
If $\lambda \in \Lambda_n^+$, then
$$\lambda p^k + \rho_{(1,\dots,1)} \cdot \frac{p^k-1}{p-1} \xmapsto{k \text{ iterations}} \underbrace{(\cdots (}_{k}\lambda) \cdots ),$$
where $\rho_{(1,\dots,1)} = (n-1, n-3, \dots, 3-n, 1-n)$.
\end{proposition}
\begin{remark}
It should be noted that $\rho_\alpha$ is more generally defined for partitions $\alpha \vdash n$, but this is not needed in this paper.
\end{remark}
\begin{proof}
Let $\rho = \rho_{(1,\dots,1)}$. One iteration of $LV_p^{-1}$ starting from $\lambda$ would return $\lambda p + \rho$, where the $\rho$ comes from the map $E$. Repeating this recursively, we obtain
$$(\cdots(((\lambda p + \rho)p + \rho)p + \rho)\cdots )p + \rho = \lambda p^k + \rho (1 + p + \cdots + p^{k-1}) = \lambda p^k + \rho \cdot \frac{p^k-1}{p-1},$$
as desired.
\end{proof}

This proposition allows us to always give one family of distinguished weights for any $n$.

\begin{corollary}\label{cor:one-sequence-distinguished}
Let $\alpha = (1, \dots, 1)$ denote a partition of $n$. Then $\rho_\alpha \cdot \frac{p^m-1}{p-1} \in \Lambda_{n, \dist}^+$ for $m \geq 0$.
\end{corollary}
\begin{proof}
We know $(0, \dots, 0) \in \Lambda_{n,\dist}^+$, then apply Proposition~\ref{prop:one_sequence_iterate}.
\end{proof}

This family of distinguished weights corresponds to the partition $(1, \dots, 1)$ of $n$. As we can already see in the $n=4$ example in Section~\ref{sec:examples}, the other families of distinguished weights can become very complicated. 

\subsection{Anti-symmetry of distinguished weights}\label{subsec:antisymmetry}

One clear observation about the distinguished weights so far is about their anti-symmetry. In particular, if $\lambda \in \Lambda_{n,\dist}^+$ with $\lambda = (\lambda_1, \dots, \lambda_n)$ then $\lambda_i = -\lambda_{n+1-i}$. We prove this statement in this section, by first defining the map $R$ for more concise notation and then proving a more general statement.

\begin{definition}
Define the map $R : \Lambda_m^+ \to \Lambda_m^+$ be defined as $R:(\lambda_1, \dots, \lambda_m) \mapsto (-\lambda_m, \dots, -\lambda_1)$ for all $m$. Overloading the notation, define the map $R : \Omega_m \to \Omega_m$ be defined as\\ $R: (\mu_1, \dots, \mu_s) \mapsto (R\mu_1, \dots, R\mu_s)$ for $m$. In other words, the map $R$ negates and reverses a weakly decreasing sequence.
\end{definition}

The fact that all distinguished weights are anti-symmetric follows directly from the following theorem, and is presented as a corollary after the theorem's proof. The proof of the following theorem requires Lemma~\ref{lem:oneclump}, which is stated and proved at the end of this section, having been moved due to length.

\begin{theorem}\label{thm:negativereverse}
If $\lambda \in \Lambda_n^+$, then $LV(R(\lambda)) = R(LV(\lambda))$.
\end{theorem}
\begin{proof}
Consider $\phi(\lambda)$ for an arbitrary $\lambda$. The elements in each row have consecutive differences with magnitude at most $1$, so the elements must begin in the same maximal clump of $\lambda$. This means that each maximal clump corresponds to some set of rows in $\phi(\lambda)$, and the rows corresponding to different maximal clumps are disjoint.

Consider the map $\phi$ on $\lambda$ and $R(\lambda)$. Let the maximal clumps of $\lambda$, ordered by largest numbers first, be $A_1, \dots, A_c$, and let $A_1$ correspond to the first $a_1$ rows of $\phi(\lambda)$, $A_2$ correspond to the next $a_2$ rows, etc. Then the maximal clumps of $R(\lambda)$, ordered by largest numbers first, are $R(A_c), \dots, R(A_1)$, with $R(A_1)$ corresponding to the last $a_1$ rows of $\phi(R(\lambda))$, $R(A_2)$ corresponding to the $a_2$ rows before that, and so on.

Consider clump $A_i$, and consider $\kappa \circ E^{-1}$ on $\phi(A_i)$ and $\phi(R(A_i))$. First, $E^{-1}$ adds constants to the elements by column, then $\kappa$ takes row sums. Lemma~\ref{lem:oneclump}, whose proof is at the end of this section due to its length, states that this theorem's statement is true for one clump, so we know that $R(LV(A_i)) = LV(R(A_i))$; then for every row in $E^{-1}(\phi(A_i))$, there is a row of the same length with the opposite row sum in $E^{-1}(\phi(R(A_i)))$ and vice versa. For each column of $\phi(A_i)$, the map $E^{-1}$ adds $(1-l, 3-l, \dots, l-1)$ to the column, where $l$ is the number of elements in the column; for the corresponding column of $\phi(R(A_i))$, the map $E^{-1}$ adds the same constant. Let $t$ be an arbitrary scalar. If, instead, $E^{-1}$ added $(1-l+t, 3-l+t, \dots, l-1+t)$ to the column in $\phi(A_i)$ and added $(1-l-t, 3-l-t, \dots, l-1-t)$ to the column in $\phi(R(A_i))$, then the property that every row in $E^{-1}(\phi(A_i))$ has a corresponding row of the same length with the opposite row sum in $E^{-1}(\phi(R(A_i)))$ (and vice versa) still holds, since all row sums are translated by $+t$ in $E^{-1}(\phi(A_i))$ and by $-t$ in $E^{-1}(\phi(R(A_i)))$.

The map $E^{-1}$ adds an anti-symmetric constant to every column, specifically in the form of $(1-L, 3-L, \dots, L-1)$, where $L$ is the total number of elements in the column. In $\phi(\lambda)$, clump $A_i$ is represented by the rows $a_1+\cdots+a_{i-1}+1$ to $a_1+\cdots+a_{i-1}+a_i$ from the top. This means that under the map $E^{-1}$, the constant that is added to the $j$th column for the rows corresponding to clump $A_i$ in $\phi(\lambda)$ is in the form $(1-l+t, 3-l+t, \dots, l-1+t)$ for some $t$, where $l$ is the number of elements in the $j$th column in clump $A_i$. On the other hand, in $\phi(R(\lambda))$, clump $R(A_i)$ is represented by the rows $a_1+\cdots+a_{i-1}+1$ to $a_1+\cdots+a_{i-1}+a_i$ from the bottom. Since the corresponding rows indices for clump $R(A_i)$ in $\phi(R(\lambda))$ are opposite to clump $A_i$ in $\phi(\lambda)$, the constant added to the $j$th column for the rows corresponding to clump $R(A_i)$ in $\phi(R(\lambda))$ is exactly in the form $(1-l-t, 3-l-t, \dots, l-1-t)$ for the same $t$. Therefore, for every row in $E^{-1}(\phi(\lambda))$, there is a corresponding row in $E^{-1}(\phi(R(\lambda)))$ of the same length with an opposite row sum and vice versa.

This means that if $\kappa(E^{-1}(\phi(\lambda))) = (\mu_1, \dots, \mu_s)$, then $\kappa(E^{-1}(\phi(R(\lambda))))$ has all row sums negated, which is exactly $(R\mu_1, \dots, R\mu_s)$. Hence, $LV(R(\lambda)) = R(LV(\lambda))$, as desired.
\end{proof}

\begin{corollary}\label{cor:antisymmetric}
If $\lambda \in \Lambda_{n,\dist}^+$, then $R(\lambda) = \lambda$, i.e.~if $\lambda = (\lambda_1, \dots, \lambda_n)$ then $\lambda_i = -\lambda_{n+1-i}$.
\end{corollary}
\begin{proof}
There exists some $m \geq 0$ such that $LV_p^m(\lambda)$ is all zeroes, which we denote simply as $LV_p^m(\lambda) = 0$. By the above theorem (Theorem~\ref{thm:negativereverse}), and since dividing by $p$ clearly commutes with $R$, we have
$$LV_p^m(\lambda) = 0 = R(0) = R(LV_p^m(\lambda)) = LV_p^m(R(\lambda)).$$
In other words, since any sequence of all zeroes is invariant under $R$, then applying $LV_p^{-1}$ preserves this. Since $LV_p$ is bijective, then $R(\lambda) = \lambda$, as desired.
\end{proof}

\begin{remark}
The converse of this theorem is false, i.e.~not all anti-symmetric sequences are distinguished weights. This can be quickly seen in the examples given in Section~\ref{sec:examples}.
\end{remark}

Finally, to finish the proof of the theorem, we must prove the following lemma. This proof includes somewhat lengthy casework, and examples of each case are given after the proof for concreteness.

Define $\phi'$, a map similar to $\phi$ with the only difference being we index the columns of the image starting at column $0$. The rest of the algorithm for constructing the map is the same, and this new indexing entails a change in parity, which leads to a slight change in the resulting image. Accordingly, define $LV' = \kappa \circ E^{-1} \circ \phi'$. This definition allows us to have a slightly stronger statement for the lemma, which is needed for induction.

\begin{lemma}\label{lem:oneclump}
Let $\lambda \in \Lambda_n^+$ be a maximal clump, i.e.~the underlying set of $\lambda$ consists of consecutive numbers. Then $R(LV(\lambda)) = LV(R(\lambda))$ and $R(LV'(\lambda)) = LV'(R(\lambda))$. 
\end{lemma}
\begin{proof}
Let the underlying set of $\lambda$ be $B$. We proceed by induction on $|B|$.

If $|B| = 1$, then $\lambda = (\lambda_1, \dots, \lambda_1)$, where $\lambda_1$ appears $m_1$ times. Then it is clear that
$$LV'(\lambda) = LV(\lambda) = ((), \dots, (), (m_1\lambda_1))$$
and
$$LV'(R(\lambda)) = LV(R(\lambda)) = ((), \dots, (), (-m_1\lambda_1)).$$
Suppose that $|B| = k$, and that the statement is true for smaller values of $|B|$. Let $B = \{\lambda_1, \dots, \lambda_k\}$, where $\lambda_i = \lambda_{i+1}+1$ for $1 \leq i \leq k-1$ and $\lambda_i$ occurs with multiplicity $m_i$ in $\lambda$.

We first show that $R(LV(\lambda)) = LV(R(\lambda))$. Since $\kappa$ takes row sums, then it suffices to show that for every row of $E^{-1}(\phi(\lambda))$, there is a corresponding row in $E^{-1}(\phi(R(\lambda)))$ with the same length but opposite row sum. Call this the row-length-sum property of $\lambda$ and $R(\lambda)$.

We split into cases:
\begin{description}
    \item[Case 1]  $k$ is odd. Let $k=2t+1$. We construct the columns of $E^{-1}(\phi(\lambda))$ one by one (instead of constructing $\phi(\lambda)$ first then applying $E^{-1}$). Each of the first $\min\{m_i : i \text{ odd}\}$ columns of $\phi(\lambda)$ are $\lambda_{\text{odd}} := (\lambda_1, \lambda_3, \dots, \lambda_{2t+1})$, and applying $E^{-1}$ makes the column have values $(\lambda_1-t, \lambda_3-(t-2), \dots, \lambda_{2t+1}+t)$. Since the $\lambda_i$'s are decreasing consecutive numbers, then
    $$s := \lambda_1-t = \lambda_3-(t-2) = \cdots = \lambda_{2t+1}+t,$$
    so the column is just copies of $s$, as shown in the diagram below.
    \vspace{-5mm}
    \[\begin{tikzcd}[cramped, column sep=tiny,row sep=small]
    	{} & {} && {} & {\phantom{m}} \\
    	&&&&&& s & s & \cdots & s \\
    	{\lambda_{\text{odd}}} & {\lambda_{\text{odd}}} & \cdots & {\lambda_{\text{odd}}} &&& \vdots & \vdots & \cdots & \vdots \\
    	&&&&&& s & s & \cdots & s \\
    	{} & {} && {}
    	\arrow[shorten >=10pt, no head, from=3-1, to=1-1]
    	\arrow[shorten >=6pt, no head, from=3-1, to=5-1]
    	\arrow[shorten >=10pt, no head, from=3-2, to=1-2]
    	\arrow[shorten >=6pt, no head, from=3-2, to=5-2]
    	\arrow[shorten >=10pt, no head, from=3-4, to=1-4]
    	\arrow["{E^{-1}}", shorten <=8pt, shorten >=8pt, maps to, from=3-4, to=3-7]
    	\arrow[shorten >=6pt, no head, from=3-4, to=5-4]
    \end{tikzcd}\]
    At this point in the algorithm of $\phi$, one of the numbers $\{\lambda_i : i \text{ odd}\}$ has no more copies left. This means that we have at least two maximal clumps now, each with underlying set size smaller than $k$. Let the new clumps be $A_1, \dots, A_c$, ordered by largest numbers first.
    
    Similarly, the first $\min\{m_i : i \text{ odd}\}$ columns of $\phi(R(\lambda))$ are\\ $R(\lambda_{\text{odd}}) := (-\lambda_{2t+1}, -\lambda_{2t-1}, \dots, -\lambda_1)$, and applying $E^{-1}$ makes the column have values $(-\lambda_{2t+1}-t, -\lambda_{2t-1}-(t-2), \dots, -\lambda_1+t)$, which are exactly copies of $-s$, as shown in the diagram below. The maximal clumps at this point in the algorithm are $R(A_c), \dots, R(A_1)$, ordered by largest numbers first.
    \vspace{-5mm}
    \[\begin{tikzcd}[cramped, column sep=tiny,row sep=small]
    	{} & {} && {} & {\phantom{m}} \\
    	&&&&&& -s & -s & \cdots & -s \\
    	{R(\lambda_{\text{odd}})} & {R(\lambda_{\text{odd}})} & \cdots & {R(\lambda_{\text{odd}})} &&& \vdots & \vdots & \cdots & \vdots \\
    	&&&&&& -s & -s & \cdots & -s \\
    	{} & {} && {}
    	\arrow[shorten >=10pt, no head, from=3-1, to=1-1]
    	\arrow[shorten >=6pt, no head, from=3-1, to=5-1]
    	\arrow[shorten >=10pt, no head, from=3-2, to=1-2]
    	\arrow[shorten >=6pt, no head, from=3-2, to=5-2]
    	\arrow[shorten >=10pt, no head, from=3-4, to=1-4]
    	\arrow["{E^{-1}}", shorten <=8pt, shorten >=8pt, maps to, from=3-4, to=3-7]
    	\arrow[shorten >=6pt, no head, from=3-4, to=5-4]
    \end{tikzcd}\]

    Let $r = \min\{m_i : i \text{ odd}\}+1$ be the next column to construct in $E^{-1}(\phi(\lambda))$. We have two further subcases:

    \item[Subcase 1a] $r$ is odd. Then $r$ is the same parity as $1$, so constructing the next column of $E^{-1}(\phi(\lambda))$ is the same as constructing the first column of $E^{-1}(\phi(A_i))$, for the clumps $A_i$, and appending them in the correct places. This is the same for $E^{-1}(\phi(R(A_i)))$. We apply the inductive hypothesis to each clump, which has underlying set size less than $k$. In particular, for each clump $A_i$, we know that for each row in $E^{-1}(\phi(A_i))$, there is a row of equal length and opposite row sum in $E^{-1}(\phi(R(A_i)))$. However, there are two differences between $E^{-1}(\phi(A_i))$ and when the same clump is constructed in $E^{-1}(\phi(\lambda))$: first, we have copies of $s$ in front in $E^{-1}(\phi(\lambda))$, and second, $E^{-1}$ may add different constants. We claim that these differences do not change the row-length-sum property of $A_i$ and $R(A_i)$:
        
    First, adding copies of $s$ in front of every row of $E^{-1}(\phi(A_i))$ and adding copies of $-s$ in front of every row of $E^{-1}(\phi(A_i))$ does not change the row-length-sum property, since row lengths are all increased by the same number and row sums are increased in $E^{-1}(\phi(A_i))$ the same amount as they are decreased in $E^{-1}(\phi(R(A_i)))$.
        
    Second, for each column $j$, $E^{-1}$ applied to $\phi(A_i)$ adds the same constants as $E^{-1}$ applied to clump $A_i$ in $\phi(\lambda)$ but shifted by a constant $c_j$. Similarly, for each column $j$, $E^{-1}$ applied to $\phi(R(A_i))$ adds the same constants as $E^{-1}$ applied to clump $R(A_i)$ in $\phi(R(\lambda))$ but shifted by $-c_j$, since the order of the clumps in $R(\lambda)$ are the reverse of those in $\lambda$. However, adding a constant to the rows in $E^{-1}(\phi(A_i))$ and subtracting the same constant to the corresponding rows in $E^{-1}(\phi(R(A_i)))$ does not change the row-length-sum property either.

    Therefore, for every clump $A_i$, the row-length-sum property is preserved when it is constructed in $E^{-1}(\phi(\lambda))$ and $E^{-1}(\phi(R(\lambda)))$. The rows in $E^{-1}(\phi(\lambda))$ which stop after column $r-1$ clearly have their corresponding row in $E^{-1}(\phi(R(\lambda)))$ which is also length $r-1$, and their sums are $(r-1)s$ and $(r-1)(-s)$, respectively.

    \item[Subcase 1b] $r$ is even. The same proof as Subcase 1a works, replacing $\phi$ with $\phi'$.

    \item[Case 2] $k$ is even. Let $k = 2t$. The first column of $\phi(\lambda)$ is $\lambda_{\text{odd}} := (\lambda_1, \lambda_3, \dots, \lambda_{2t-1})$. The second column is $\lambda_{\text{even}} := (\lambda_2, \lambda_4, \dots, \lambda_{2t})$. The third column is $\lambda_{\text{odd}}$, and in general, the columns alternate until some $\lambda_i$ has no more copies left.
    
    Similarly, the first column in $\phi(R(\lambda))$ is $(-\lambda_{2t}, -\lambda_{2t-2}, \dots, -\lambda_2) = R(\lambda_{\text{even}})$; the second column is $(-\lambda_{2t-1}, -\lambda_{2t-3}, \dots, -\lambda_1) = R(\lambda_{\text{odd}})$; the third column is $R(\lambda_{\text{even}})$. The columns alternate until some number has no copies left.
    
    Once some $\lambda_i$ has no more copies left in $\phi(\lambda)$, then to construct the next column, the remaining numbers split into more than one maximal clump. Let the $(r-1)$th column be the column after which, the remaining numbers are no longer one clump. We split into the following subcases:

    \item[Subcase 2a] $r$ is even. This means that $r-1$ is odd, so there are no more copies of some $\lambda_i$ with $i$ odd. Let $C$ be the set of these indices which have no more copies, i.e.
    $$C = \{i : i \text{ odd}, m_i = \frac r2\} = \{c_1 \leq \cdots \leq c_d\}.$$
    We claim that the $r$th column of $\phi(\lambda)$ is $\lambda_{\text{even}}$. The underlying sets of the maximal clumps are
    $$\{\lambda_1, \dots, \lambda_{c_1-1}\}, \{\lambda_{c_1+1}, \dots, \lambda_{c_2-1}\}, \dots, \{\lambda_{c_d+1}, \dots, \lambda_{2t}\}.$$
    The first clump's underlying set has an even size. Since $r$ is even, the selected elements are $\lambda_2, \lambda_4, \dots, \lambda_{c_1-1}$. The $i$th clump, for $2 \leq i \leq d$, has underlying set $\{\lambda_{c_{i-1}+1}, \dots, \lambda_{c_i-1}\}$, which has odd size, so the selected elements are $\lambda_{c_{i-1}+1}, \lambda_{c_{i-1}+3}, \dots, \lambda_{c_i-1}$. The $(d+1)$th clump has underlying set $\{\lambda_{c_d+1}, \dots, \lambda_{2t}\}$, which has odd size, so the selected elements are $\lambda_{c_d+1}, \lambda_{c_d+3}, \dots, \lambda_{2t}$. Therefore, to construct the $r$th column, all of the even indices are selected, so the $r$th column is $(\lambda_2, \lambda_4, \dots, \lambda_{2t}) = \lambda_{\text{even}}$.

    Therefore, the first $r$ columns of $\phi(\lambda)$ are $\lambda_{\text{odd}}, \lambda_{\text{even}}, \lambda_{\text{odd}}, \dots, \lambda_{\text{even}}$. Applying $E^{-1}$ makes $\lambda_{\text{odd}}$ to have values $(\lambda_1-(t-1), \lambda_3-(t-3), \dots, \lambda_{2t-1}+(t-1))$ and $\lambda_{\text{even}}$ to have values $(\lambda_2-(t-1), \lambda_4-(t-3), \dots, \lambda_{2t}+(t-1))$. Since the $\lambda_i$'s are decreasing consecutive numbers, we have that
    $$s := \lambda_1-(t-1) = \lambda_3 - (t-3) = \cdots = \lambda_{2t-1}+(t-1)$$
    and
    $$s-1 = \lambda_2-(t-1) = \lambda_4 - (t-3) = \cdots = \lambda_{2t} + (t-1).$$
    Hence, the columns alternate between all $s$'s and all $(s-1)$'s, as shown in the diagram below.
    \vspace{-5mm}
    \[\begin{tikzcd}[sep=tiny]
    	{} & {} & {} && {} & {\phantom{m}} \\
    	&&&&&&& s & {s-1} & s & \cdots & {s-1} \\
    	{\lambda_{\text{odd}}} & {\lambda_{\text{even}}} & {\lambda_{\text{odd}}} & \cdots & {\lambda_{\text{even}}} &&& \vdots & \vdots & \vdots & \cdots & \vdots \\
    	&&&&&&& s & {s-1} & s & \cdots & {s-1} \\
    	{} & {} & {} && {}
    	\arrow[shorten >=10pt, no head, from=3-1, to=1-1]
    	\arrow[shorten >=10pt, no head, from=3-1, to=5-1]
    	\arrow[shorten >=10pt, no head, from=3-2, to=1-2]
    	\arrow[shorten >=6pt, no head, from=3-2, to=5-2]
    	\arrow[shorten >=10pt, no head, from=3-3, to=1-3]
    	\arrow[shorten >=6pt, no head, from=3-3, to=5-3]
    	\arrow[shorten >=10pt, no head, from=3-5, to=1-5]
    	\arrow["{{E^{-1}}}", shorten <=8pt, shorten >=8pt, maps to, from=3-5, to=3-8]
    	\arrow[shorten >=6pt, no head, from=3-5, to=5-5]
    \end{tikzcd}\]
    Now, we examine what occurs in $\phi(R(\lambda))$. The $(r-2)$th column is $R(\lambda_{\text{odd}})$, and by assumption, there are still copies of all $\lambda_i$ remaining. The $(r-1)$th column is $R(\lambda_{\text{even}})$. If there are still copies of all $\lambda_i$ remaining, then the $r$th column is $R(\lambda_{\text{odd}})$. Otherwise, after the $(r-1)$th column, at least one even index $i$ has no more copies of $\lambda_i$. Let $C$ be the set of these indices which have no more copies, i.e. $C = \{i : i \text{ even}, m_i = \frac{r}2\} = \{c_1\geq \cdots \geq c_d\}$. We claim that the $r$th column is still $R(\lambda_{\text{odd}})$. The underlying sets of the maximal clumps are
    $$\{-\lambda_{2t}, \dots, -\lambda_{c_d+1}\}, \{-\lambda_{c_d-1}, \dots, -\lambda_{c_{d-1}+1}\}, \dots, \{-\lambda_{c_1-1}, \dots, -\lambda_1\}.$$
    The first clump's underlying set has an even size. Since $r$ is even, the selected elements are $-\lambda_{2t-1}, -\lambda_{2t-3}, \dots, -\lambda_{c_1+1}$. The $i$th clump, for $2 \leq i \leq d$, has underlying set\\ $\{-\lambda_{c_{d-i+2}-1}, \dots, -\lambda_{c_{d-i+1}+1}\}$, which has odd size, so the selected elements are exactly $-\lambda_{c_{d-i+2}-1}, -\lambda_{c_{d-i+2}-3}, \dots, -\lambda_{c_{d-i+1}+1}$. The $(d+1)$th clump has underlying set\\ $\{-\lambda_{c_d-1}, \dots, -\lambda_1\}$, which has odd size, so the selected elements are $-\lambda_{c_d-1}, -\lambda_{c_d-3}, \dots, -\lambda_1$. Therefore, to construct the $r$th column, all of the odd indices are selected, so the $r$th column is still $(-\lambda_{2t-1}, -\lambda_{2t-3}, \dots, -\lambda_1) = R(\lambda_{\text{odd}})$.

    In both cases, the $r$th column of $\phi(R(\lambda))$ is $R(\lambda_{\text{odd}})$. Therefore, the first $r$ columns of $\phi(R(\lambda))$ are $R(\lambda_{\text{even}}), R(\lambda_{\text{odd}}), R(\lambda_{\text{even}}), \dots, R(\lambda_{\text{odd}})$. Applying $E^{-1}$ makes $R(\lambda_{\text{even}})$ to have values $(-\lambda_{2t}-(t-1), -\lambda_{2t-2}-(t-3), \dots, -\lambda_2+(t-1))$ and $R(\lambda_{\text{odd}})$ to have values $(-\lambda_{2t-1}-(t-1), -\lambda_{2t-3}-(t-3), \dots, -\lambda_1+(t-1))$, which satisfies
    $$-(s-1) = -\lambda_{2t}-(t-1) = -\lambda_{2t-2}-(t-3) = \cdots = -\lambda_2 + (t-1)$$
    and
    $$-s = -\lambda_{2t-1}-(t-1) = -\lambda_{2t-3}-(t-3) = \cdots = -\lambda_1 + (t-1).$$
    Hence, the columns alternate between all $-(s-1)$'s and all $-s$'s, as shown in the diagram below.
    \vspace{-5mm}
    {\small
    \[\begin{tikzcd}[sep=tiny]
    	{} & {} & {} && {} & {\phantom{m}} \\
    	&&&&&&& -(s-1) & {-s} & -(s-1) & \cdots & {-s} \\
    	{R(\lambda_{\text{even}})} & {R(\lambda_{\text{odd}})} & {R(\lambda_{\text{even}})} & \cdots & {R(\lambda_{\text{odd}})} &&& \vdots & \vdots & \vdots & \cdots & \vdots \\
    	&&&&&&& -(s-1) & {-s} & -(s-1) & \cdots & {-s} \\
    	{} & {} & {} && {}
    	\arrow[shorten >=10pt, no head, from=3-1, to=1-1]
    	\arrow[shorten >=6pt, no head, from=3-1, to=5-1]
    	\arrow[shorten >=10pt, no head, from=3-2, to=1-2]
    	\arrow[shorten >=6pt, no head, from=3-2, to=5-2]
    	\arrow[shorten >=10pt, no head, from=3-3, to=1-3]
    	\arrow[shorten >=6pt, no head, from=3-3, to=5-3]
    	\arrow[shorten >=10pt, no head, from=3-5, to=1-5]
    	\arrow["{{E^{-1}}}", shorten <=8pt, shorten >=8pt, maps to, from=3-5, to=3-8]
    	\arrow[shorten >=6pt, no head, from=3-5, to=5-5]
    \end{tikzcd}\]
    }
    Therefore, the first $r$ columns of $E^{-1}(\phi(\lambda))$ alternate between all $s$'s and all $s-1$'s, while the first $r$ columns of $E^{-1}(\phi(R(\lambda)))$ alternate between all $-(s-1)$'s and all $-s$'s. In any case, the elements in a given column are the same, and since the clumps to construct the $(r+1)$th column have underlying set size less than $k$, we apply the inductive hypothesis in the same way as Subcase 1a, since $r+1$ is odd.

    \item[Subcase 2b] $r$ is odd. This means that $r-1$ is even, so there are no more copies of some $\lambda_i$ with $i$ even. Let $C$ be the set of these indices which have no more copies, i.e.
    $$C = \{i : i \text{ even}, m_i = \frac{r-1}2\} = \{c_1 \geq \cdots \geq c_d\}.$$
    Consider what occurs in $\phi(R(\lambda))$. The $(r-2)$th column is $R(\lambda_{\text{even}})$, and by assumption, there are some $\lambda_i$ with no more copies remaining, so we have more than one clump when constructing column $r-1$. We claim that the $(r-1)$th column is $R(\lambda_{\text{odd}})$. The underlying sets of the maximal clumps are
    $$\{-\lambda_{2t}, \dots, -\lambda_{c_d+1}\}, \{-\lambda_{c_d-1}, \dots, -\lambda_{c_{d-1}+1}\}, \dots, \{-\lambda_{c_1-1}, \dots, -\lambda_1\}.$$
    The first clump's underlying set has even size. Since $r-1$ is even, the selected elements are exactly $-\lambda_{2t-1}, -\lambda_{2t-3}, \dots, -\lambda_{c_1+1}$. The $i$th clump, for $2 \leq i \leq d$, has underlying set $\{-\lambda_{c_{d-i+2}-1}, \dots, -\lambda_{c_{d-i+1}+1}\}$, which has odd size, so the selected elements are $-\lambda_{c_{i-1}-1}, -\lambda_{c_{i-1}-3}, \dots, -\lambda_{c_i+1}$. The $(d+1)$th clump has underlying set $\{-\lambda_{c_d-1}, \dots, -\lambda_1\}$, which has odd size, so the selected elements are $-\lambda_{c_d-1}, -\lambda_{c_d-3}, \dots, -\lambda_1$. Therefore, in construction of the $(r-1)$th column, all of the odd indices are selected, so the $(r-1)$th column is $(-\lambda_{2t-1}, -\lambda_{2t-3}, \dots, -\lambda_1) = R(\lambda_{\text{odd}})$, as desired.

    Therefore, the first $r-1$ columns of $\phi(\lambda)$ are $\lambda_{\text{odd}}, \lambda_{\text{even}}, \lambda_{\text{odd}}, \dots, \lambda_{\text{even}}$, and the first $r-1$ columns of $\phi(R(\lambda))$ are $R(\lambda_{\text{even}}), R(\lambda_{\text{odd}}), R(\lambda_{\text{even}}), \dots, R(\lambda_{\text{odd}})$. As in Subcase 2a, the first $r-1$ columns of $E^{-1}(\phi(\lambda))$ alternates between all $s$'s and all $(s-1)$'s, and the first $r-1$ columns of $E^{-1}(\phi(R(\lambda)))$ alternates between all $-(s-1)$'s and all $-s$'s, which is enough to apply the inductive hypothesis in the same way as Subcase 1a, since all clumps have underlying set of size less than $k$ and $r$ is odd.
\end{description}

Finally, we want to show $R(LV'(\lambda)) = LV'(R(\lambda))$. This is essentially the same casework as above, replacing $\phi$ with $\phi'$. By induction, we are done.
\end{proof}

\begin{example}
Here are examples for all of the subcases that appear above. It should be noted that the numbers themselves do not affect the shape of the image, but rather the multiplicities of the numbers do.
\begin{description}
    \newpage
    \item[Subcase 1a] Let $\lambda = (9,9,9,8,8,7,7,6,6,5,5,5,5,4,4,4,3,3)$. Then $r=3$, and we have
    {\small
    \begin{center}
        \begin{tikzpicture}[scale=0.55]
            \node at (-4,0) {$\lambda \xmapsto{\phi}$};
            \draw (-3,1) rectangle node{$9$} ++(1,1);
            \draw (-3,0) rectangle node{$7$} ++(1,1);
            \draw (-3,-1) rectangle node{$5$} ++(1,1);
            \draw (-3,-2) rectangle node{$3$} ++(1,1);
            \draw (-2,1) rectangle node{$9$} ++(1,1);
            \draw (-2,0) rectangle node{$7$} ++(1,1);
            \draw (-2,-1) rectangle node{$5$} ++(1,1);
            \draw (-2,-2) rectangle node{$3$} ++(1,1);
            \draw (-1,1) rectangle node{$9$} ++(1,1);
            \draw (-1,-1) rectangle node{$6$} ++(1,1);
            \draw (-1,-2) rectangle node{$4$} ++(1,1);
            \draw (0,1) rectangle node{$8$} ++(1,1);
            \draw (0,-1) rectangle node{$6$} ++(1,1);
            \draw (0,-2) rectangle node{$4$} ++(1,1);
            \draw (1,1) rectangle node{$8$} ++(1,1);
            \draw (1,-2) rectangle node{$5$} ++(1,1);
            \draw (2,-2) rectangle node{$4$} ++(1,1);
            \draw (3,-2) rectangle node{$5$} ++(1,1);
            \node at (-2,-2.5) {$\underbrace{\phantom{mmm}}_{\min\{m_i: i\text{ odd}\}}$};
            \node at (-0.5,3) {$\overset{r}{\downarrow}$};
            \node at (5,0) {$\xmapsto{E^{-1}}$};
            \draw (6,1) rectangle node{$6$} ++(1,1);
            \draw (6,0) rectangle node{$6$} ++(1,1);
            \draw (6,-1) rectangle node{$6$} ++(1,1);
            \draw (6,-2) rectangle node{$6$} ++(1,1);
            \draw (7,1) rectangle node{$6$} ++(1,1);
            \draw (7,0) rectangle node{$6$} ++(1,1);
            \draw (7,-1) rectangle node{$6$} ++(1,1);
            \draw (7,-2) rectangle node{$6$} ++(1,1);
            \draw (8,1) rectangle node{$7$} ++(1,1);
            \draw (8,-1) rectangle node{$6$} ++(1,1);
            \draw (8,-2) rectangle node{$6$} ++(1,1);
            \draw (9,1) rectangle node{$6$} ++(1,1);
            \draw (9,-1) rectangle node{$6$} ++(1,1);
            \draw (9,-2) rectangle node{$6$} ++(1,1);
            \draw (10,1) rectangle node{$7$} ++(1,1);
            \draw (10,-2) rectangle node{$6$} ++(1,1);
            \draw (11,-2) rectangle node{$4$} ++(1,1);
            \draw (12,-2) rectangle node{$5$} ++(1,1);
            \node at (18,0) {$\xmapsto{\kappa} ((), (12), (), (24), (32), (), (39))$};
        \end{tikzpicture}
    \end{center}
    \medskip
    \begin{center}
        \begin{tikzpicture}[scale=0.55]
            \node at (-4.5,0) {$R(\lambda) \xmapsto{\phi}$};
            \draw (-3,1) rectangle node{$-3$} ++(1,1);
            \draw (-3,0) rectangle node{$-5$} ++(1,1);
            \draw (-3,-1) rectangle node{$-7$} ++(1,1);
            \draw (-3,-2) rectangle node{$-9$} ++(1,1);
            \draw (-2,1) rectangle node{$-3$} ++(1,1);
            \draw (-2,0) rectangle node{$-5$} ++(1,1);
            \draw (-2,-1) rectangle node{$-7$} ++(1,1);
            \draw (-2,-2) rectangle node{$-9$} ++(1,1);
            \draw (-1,0) rectangle node{$-4$} ++(1,1);
            \draw (-1,-1) rectangle node{$-6$} ++(1,1);
            \draw (-1,-2) rectangle node{$-8$} ++(1,1);
            \draw (0,0) rectangle node{$-4$} ++(1,1);
            \draw (0,-1) rectangle node{$-6$} ++(1,1);
            \draw (0,-2) rectangle node{$-9$} ++(1,1);
            \draw (1,0) rectangle node{$-4$} ++(1,1);
            \draw (1,-2) rectangle node{$-8$} ++(1,1);
            \draw (2,0) rectangle node{$-5$} ++(1,1);
            \draw (3,0) rectangle node{$-5$} ++(1,1);
            \node at (5,0) {$\xmapsto{E^{-1}}$};
            \draw (6,1) rectangle node{$-6$} ++(1,1);
            \draw (6,0) rectangle node{$-6$} ++(1,1);
            \draw (6,-1) rectangle node{$-6$} ++(1,1);
            \draw (6,-2) rectangle node{$-6$} ++(1,1);
            \draw (7,1) rectangle node{$-6$} ++(1,1);
            \draw (7,0) rectangle node{$-6$} ++(1,1);
            \draw (7,-1) rectangle node{$-6$} ++(1,1);
            \draw (7,-2) rectangle node{$-6$} ++(1,1);
            \draw (8,0) rectangle node{$-6$} ++(1,1);
            \draw (8,-1) rectangle node{$-6$} ++(1,1);
            \draw (8,-2) rectangle node{$-6$} ++(1,1);
            \draw (9,0) rectangle node{$-6$} ++(1,1);
            \draw (9,-1) rectangle node{$-6$} ++(1,1);
            \draw (9,-2) rectangle node{$-7$} ++(1,1);
            \draw (10,0) rectangle node{$-5$} ++(1,1);
            \draw (10,-2) rectangle node{$-7$} ++(1,1);
            \draw (11,0) rectangle node{$-5$} ++(1,1);
            \draw (12,0) rectangle node{$-5$} ++(1,1);
            \node at (16.9,0) {$\xmapsto{\kappa} ((), (-12), (), (-24),$};
            \node at (18.4,-1) {$(-32), (), (-39))$};
        \end{tikzpicture}
    \end{center}
    }

    \item[Subcase 1b] Let $\lambda = (9,8,8,8,7,7,6,6,5,4,3,3)$. Then $r=2$, and we have
    {\small
    \begin{center}
        \begin{tikzpicture}[scale=0.55]
            \node at (-5,0) {$\lambda \xmapsto{\phi}$};
            \draw (-3,1) rectangle node{$9$} ++(1,1);
            \draw (-3,0) rectangle node{$7$} ++(1,1);
            \draw (-3,-1) rectangle node{$5$} ++(1,1);
            \draw (-3,-2) rectangle node{$3$} ++(1,1);
            \draw (-2,1) rectangle node{$8$} ++(1,1);
            \draw (-2,0) rectangle node{$6$} ++(1,1);
            \draw (-2,-2) rectangle node{$3$} ++(1,1);
            \draw (-1,1) rectangle node{$8$} ++(1,1);
            \draw (-1,0) rectangle node{$6$} ++(1,1);
            \draw (-1,-2) rectangle node{$4$} ++(1,1);
            \draw (0,1) rectangle node{$7$} ++(1,1);
            \draw (1,1) rectangle node{$8$} ++(1,1);
            \node at (-2.5,-2.5) {$\underbrace{\phantom{}}_{\min\{m_i: i\text{ odd}\}}$};
            \node at (-1.5,3) {$\overset{r}{\downarrow}$};
            \node at (3,0) {$\xmapsto{E^{-1}}$};
            \draw (4,1) rectangle node{$6$} ++(1,1);
            \draw (4,0) rectangle node{$6$} ++(1,1);
            \draw (4,-1) rectangle node{$6$} ++(1,1);
            \draw (4,-2) rectangle node{$6$} ++(1,1);
            \draw (5,1) rectangle node{$6$} ++(1,1);
            \draw (5,0) rectangle node{$6$} ++(1,1);
            \draw (5,-2) rectangle node{$5$} ++(1,1);
            \draw (6,1) rectangle node{$6$} ++(1,1);
            \draw (6,0) rectangle node{$6$} ++(1,1);
            \draw (6,-2) rectangle node{$6$} ++(1,1);
            \draw (7,1) rectangle node{$7$} ++(1,1);
            \draw (8,1) rectangle node{$8$} ++(1,1);
            \node at (13.5,0) {$\xmapsto{\kappa} ((6), (), (18, 17), (), (33))$};
        \end{tikzpicture}
    \end{center}
    \medskip
    \begin{center}
        \begin{tikzpicture}[scale=0.55]
            \node at (-4.5,0) {$R(\lambda) \xmapsto{\phi}$};
            \draw (-3,1) rectangle node{$-3$} ++(1,1);
            \draw (-3,0) rectangle node{$-5$} ++(1,1);
            \draw (-3,-1) rectangle node{$-7$} ++(1,1);
            \draw (-3,-2) rectangle node{$-9$} ++(1,1);
            \draw (-2,1) rectangle node{$-4$} ++(1,1);
            \draw (-2,0) rectangle node{$-6$} ++(1,1);
            \draw (-2,-1) rectangle node{$-8$} ++(1,1);
            \draw (-1,1) rectangle node{$-3$} ++(1,1);
            \draw (-1,0) rectangle node{$-6$} ++(1,1);
            \draw (-1,-1) rectangle node{$-8$} ++(1,1);
            \draw (0,-1) rectangle node{$-8$} ++(1,1);
            \draw (1,-1) rectangle node{$-7$} ++(1,1);
            \node at (3,0) {$\xmapsto{E^{-1}}$};
            \draw (4,1) rectangle node{$-6$} ++(1,1);
            \draw (4,0) rectangle node{$-6$} ++(1,1);
            \draw (4,-1) rectangle node{$-6$} ++(1,1);
            \draw (4,-2) rectangle node{$-6$} ++(1,1);
            \draw (5,1) rectangle node{$-6$} ++(1,1);
            \draw (5,0) rectangle node{$-6$} ++(1,1);
            \draw (5,-1) rectangle node{$-6$} ++(1,1);
            \draw (6,1) rectangle node{$-5$} ++(1,1);
            \draw (6,0) rectangle node{$-6$} ++(1,1);
            \draw (6,-1) rectangle node{$-6$} ++(1,1);
            \draw (7,-1) rectangle node{$-8$} ++(1,1);
            \draw (8,-1) rectangle node{$-7$} ++(1,1);
            \node at (14,0) {$\xmapsto{\kappa} ((-6), (), (-17, -18), (), (-33))$};
        \end{tikzpicture}
    \end{center}
    }

    \item[Subcase 2a] Let $\lambda = (9,9,8,8,7,6,6,5,5,4,4,4)$. Then $r=2$, and we have
    {\small
    \begin{center}
        \begin{tikzpicture}[scale=0.55]
            \node at (-5,0.5) {$\lambda \xmapsto{\phi}$};
            \draw (-3,1) rectangle node{$9$} ++(1,1);
            \draw (-3,0) rectangle node{$7$} ++(1,1);
            \draw (-3,-1) rectangle node{$5$} ++(1,1);
            \draw (-2,1) rectangle node{$8$} ++(1,1);
            \draw (-2,0) rectangle node{$6$} ++(1,1);
            \draw (-2,-1) rectangle node{$4$} ++(1,1);
            \draw (-1,1) rectangle node{$9$} ++(1,1);
            \draw (-1,0) rectangle node{$6$} ++(1,1);
            \draw (-1,-1) rectangle node{$4$} ++(1,1);
            \draw (0,1) rectangle node{$8$} ++(1,1);
            \draw (0,-1) rectangle node{$4$} ++(1,1);
            \draw (1,-1) rectangle node{$5$} ++(1,1);
            \node at (-1.5,3) {$\overset{r}{\downarrow}$};
            \node at (3,0.5) {$\xmapsto{E^{-1}}$};
            \draw (4,1) rectangle node{$7$} ++(1,1);
            \draw (4,0) rectangle node{$7$} ++(1,1);
            \draw (4,-1) rectangle node{$7$} ++(1,1);
            \draw (5,1) rectangle node{$6$} ++(1,1);
            \draw (5,0) rectangle node{$6$} ++(1,1);
            \draw (5,-1) rectangle node{$6$} ++(1,1);
            \draw (6,1) rectangle node{$7$} ++(1,1);
            \draw (6,0) rectangle node{$6$} ++(1,1);
            \draw (6,-1) rectangle node{$6$} ++(1,1);
            \draw (7,1) rectangle node{$7$} ++(1,1);
            \draw (7,-1) rectangle node{$5$} ++(1,1);
            \draw (8,-1) rectangle node{$5$} ++(1,1);
            \node at (13,0.5) {$\xmapsto{\kappa} ((), (), (19), (27), (29))$};
        \end{tikzpicture}
    \end{center}
    \medskip
    \begin{center}
        \begin{tikzpicture}[scale=0.55]
            \node at (-5,0.5) {$R(\lambda) \xmapsto{\phi}$};
            \draw (-3,1) rectangle node{$-4$} ++(1,1);
            \draw (-3,0) rectangle node{$-6$} ++(1,1);
            \draw (-3,-1) rectangle node{$-8$} ++(1,1);
            \draw (-2,1) rectangle node{$-5$} ++(1,1);
            \draw (-2,0) rectangle node{$-7$} ++(1,1);
            \draw (-2,-1) rectangle node{$-9$} ++(1,1);
            \draw (-1,1) rectangle node{$-4$} ++(1,1);
            \draw (-1,0) rectangle node{$-6$} ++(1,1);
            \draw (-1,-1) rectangle node{$-8$} ++(1,1);
            \draw (0,1) rectangle node{$-5$} ++(1,1);
            \draw (0,-1) rectangle node{$-9$} ++(1,1);
            \draw (1,1) rectangle node{$-4$} ++(1,1);
            \node at (-1.5,3) {$\overset{r}{\downarrow}$};
            \node at (-2.5,3) {$\overset{r-1}{\downarrow}$};
            \node at (3,0.5) {$\xmapsto{E^{-1}}$};
            \draw (4,1) rectangle node{$-6$} ++(1,1);
            \draw (4,0) rectangle node{$-6$} ++(1,1);
            \draw (4,-1) rectangle node{$-6$} ++(1,1);
            \draw (5,1) rectangle node{$-7$} ++(1,1);
            \draw (5,0) rectangle node{$-7$} ++(1,1);
            \draw (5,-1) rectangle node{$-7$} ++(1,1);
            \draw (6,1) rectangle node{$-6$} ++(1,1);
            \draw (6,0) rectangle node{$-6$} ++(1,1);
            \draw (6,-1) rectangle node{$-6$} ++(1,1);
            \draw (7,1) rectangle node{$-6$} ++(1,1);
            \draw (7,-1) rectangle node{$-8$} ++(1,1);
            \draw (8,1) rectangle node{$-4$} ++(1,1);
            \node at (14,0.5) {$\xmapsto{\kappa} ((), (), (-19), (-27), (-29))$};
        \end{tikzpicture}
    \end{center}
    }

    \newpage
    \item[Subcase 2b] Let $\lambda = (9,9,8,8,7,7,7,6,5,5,4,4)$. Then $r=3$, and we have
    {\small
    \begin{center}
        \begin{tikzpicture}[scale=0.55]
            \node at (-5,0.5) {$\lambda \xmapsto{\phi}$};
            \draw (-3,1) rectangle node{$9$} ++(1,1);
            \draw (-3,0) rectangle node{$7$} ++(1,1);
            \draw (-3,-1) rectangle node{$5$} ++(1,1);
            \draw (-2,1) rectangle node{$8$} ++(1,1);
            \draw (-2,0) rectangle node{$6$} ++(1,1);
            \draw (-2,-1) rectangle node{$4$} ++(1,1);
            \draw (-1,1) rectangle node{$9$} ++(1,1);
            \draw (-1,0) rectangle node{$7$} ++(1,1);
            \draw (-1,-1) rectangle node{$5$} ++(1,1);
            \draw (0,0) rectangle node{$7$} ++(1,1);
            \draw (0,-1) rectangle node{$4$} ++(1,1);
            \draw (1,0) rectangle node{$8$} ++(1,1);
            \node at (-0.5,3) {$\overset{r}{\downarrow}$};
            \node at (3,0.5) {$\xmapsto{E^{-1}}$};
            \draw (4,1) rectangle node{$7$} ++(1,1);
            \draw (4,0) rectangle node{$7$} ++(1,1);
            \draw (4,-1) rectangle node{$7$} ++(1,1);
            \draw (5,1) rectangle node{$6$} ++(1,1);
            \draw (5,0) rectangle node{$6$} ++(1,1);
            \draw (5,-1) rectangle node{$6$} ++(1,1);
            \draw (6,1) rectangle node{$7$} ++(1,1);
            \draw (6,0) rectangle node{$7$} ++(1,1);
            \draw (6,-1) rectangle node{$7$} ++(1,1);
            \draw (7,0) rectangle node{$6$} ++(1,1);
            \draw (7,-1) rectangle node{$5$} ++(1,1);
            \draw (8,0) rectangle node{$8$} ++(1,1);
            \node at (13,0.5) {$\xmapsto{\kappa} ((), (), (20), (25), (34))$};
        \end{tikzpicture}
    \end{center}
    \medskip
    \begin{center}
        \begin{tikzpicture}[scale=0.55]
            \node at (-5,0.5) {$R(\lambda) \xmapsto{\phi}$};
            \draw (-3,1) rectangle node{$-4$} ++(1,1);
            \draw (-3,0) rectangle node{$-6$} ++(1,1);
            \draw (-3,-1) rectangle node{$-8$} ++(1,1);
            \draw (-2,1) rectangle node{$-5$} ++(1,1);
            \draw (-2,0) rectangle node{$-7$} ++(1,1);
            \draw (-2,-1) rectangle node{$-9$} ++(1,1);
            \draw (-1,1) rectangle node{$-4$} ++(1,1);
            \draw (-1,0) rectangle node{$-7$} ++(1,1);
            \draw (-1,-1) rectangle node{$-9$} ++(1,1);
            \draw (0,1) rectangle node{$-5$} ++(1,1);
            \draw (0,0) rectangle node{$-8$} ++(1,1);
            \draw (1,0) rectangle node{$-7$} ++(1,1);
            \node at (-0.5,3) {$\overset{r}{\downarrow}$};
            \node at (-1.5,3) {$\overset{r-1}{\downarrow}$};
            \node at (3,0.5) {$\xmapsto{E^{-1}}$};
            \draw (4,1) rectangle node{$-6$} ++(1,1);
            \draw (4,0) rectangle node{$-6$} ++(1,1);
            \draw (4,-1) rectangle node{$-6$} ++(1,1);
            \draw (5,1) rectangle node{$-7$} ++(1,1);
            \draw (5,0) rectangle node{$-7$} ++(1,1);
            \draw (5,-1) rectangle node{$-7$} ++(1,1);
            \draw (6,1) rectangle node{$-6$} ++(1,1);
            \draw (6,0) rectangle node{$-7$} ++(1,1);
            \draw (6,-1) rectangle node{$-7$} ++(1,1);
            \draw (7,1) rectangle node{$-6$} ++(1,1);
            \draw (7,0) rectangle node{$-7$} ++(1,1);
            \draw (8,0) rectangle node{$-7$} ++(1,1);
            \node at (14,0.5) {$\xmapsto{\kappa} ((), (), (-20), (-25), (-34))$};
        \end{tikzpicture}
    \end{center}
    }
\end{description}
\end{example}

This anti-symmetry property implies that each distinguished weight is characterized by the first $\left\lfloor \frac n2 \right\rfloor$ elements. This is important when we find the asymptotics regarding $|\Lambda_{n,k}^+|$ next.

\subsection{Sizes of $\Lambda_{n,k}^+$}\label{subsec:asymptotics}

Recall that $\Lambda_{n,k}^+$ is the subset of distinguished weights that go to all zeroes with at most $k$ iterations of the algorithm.

\begin{theorem}\label{thm:lambda_nk_recursion}
We have the following recursive formula
$$|\Lambda_{n,k}^+| = 1 + k + \sum_{\substack{\alpha \vdash n \\ \alpha = (a^{\ell_a},\dots, 1^{\ell_1}) \\ \alpha \neq (n), (1, \dots, 1)}} \sum_{m=0}^{k-1} \prod_{i=1}^{a} |\Lambda_{\ell_i, m}^+|,$$
where $|\Lambda_{0,k}^+| = 1$ and $|\Lambda_{1,k}^+| = 1$ for $k \geq 0$. 
\end{theorem}
\begin{proof}
Suppose we had a sequence $\lambda \in \Lambda_{n,k}^+$. In the sequences of outputs of repeated iterations of $LV_p$, there exists $m$ such that after $k-m$ iterations it is in the form
$$\underbrace{\Big( \cdots \Big(}_{k-m} \underbrace{(*, \dots, *)}_{\ell_1}, \underbrace{(*, \dots, *)}_{\ell_2}, \dots, \underbrace{(*, \dots, *)}_{\ell_a} \Big) \cdots \Big),$$
where the $*$ are integers and there exists $\alpha \vdash n$ such that $\alpha = (a^{\ell_a},\dots,1^{\ell_1})$. If $\alpha \neq (n), (1, \dots, 1)$, then the number of such sequences, for a fixed $m$, is $\prod_{i=1}^a |\Lambda_{\ell_i, m}^+|$, because recursively there are $|\Lambda_{\ell_i, m}^+|$ sequences for $\underbrace{(*, \dots, *)}_{\ell_i}$ to go to all zeroes after at most $m$ iterations of the algorithm.

In the case of $\alpha = (1, \dots, 1)$, this $m$ is not unique, but we know that the number of sequences that go to all zeroes after at most $k$ iterations is $k+1$, given by
$$\underbrace{(\cdots(}_{m} 0, \dots, 0) \cdots )$$
for $1 \leq m \leq k+1$. In the case of $\alpha = (n)$, we know that $LV_p$ sends $(0, \dots, 0) \mapsto ((), \dots, (), (0))$, so this does not generate any additional sequences.

Therefore, the count of $|\Lambda_{n,k}^+|$ is equal to the sum
$$(k+1) + \sum_{\substack{\alpha \vdash n \\ \alpha = (a^{\ell_a},\dots,1^{\ell_1}) \\ \alpha \neq (n), (1, \dots, 1)}} \sum_{m=0}^{k-1} \prod_{i=1}^{a} |\Lambda_{\ell_i, m}^+|,$$
as desired.
\end{proof}

\begin{example}
Using this recursive formula, we can compute that
\begin{align*}
    |\Lambda_{1,k}^+| &= 1, \\
    |\Lambda_{2,k}^+| &= k+1, \\
    |\Lambda_{3,k}^+| &= 2k+1, \\
    |\Lambda_{4,k}^+| &= k^2 + 3k + 1, \\
    |\Lambda_{5,k}^+| &= 2k^2+4k+1, \\
    |\Lambda_{6,k}^+| &= \frac{4k^3 + 27k^2 + 29k + 6}6,
\end{align*}
which aligns with our explicit characterizations of the $n\leq 4$ cases.
\end{example}

From these smaller examples, we can observe that the leading power is always $\left\lfloor \frac n2 \right\rfloor$.

\begin{proposition}\label{prop:lambda_nk_asymptotic}
For a fixed $n$, we have
$$|\Lambda_{n,k}^+| = \Theta(k^{\left\lfloor \frac n2 \right\rfloor}).$$
\end{proposition}
\begin{proof}
We proceed by induction. The cases $k=1, \dots, 6$ are shown in the above example. Now, suppose that it is true for $k=1, \dots, j-1$.

Let $\alpha \vdash n$ such that $\alpha \neq (n), (1, \dots, 1)$ and $\alpha = (a^{\ell_a},\dots,1^{\ell_1})$. Then by casework on the parity of $n$, it is clear that
$$\left\lfloor \frac{\ell_1}2 \right\rfloor + \cdots + \left\lfloor \frac{\ell_a}2 \right\rfloor \leq \left\lfloor \frac{n-2}2 \right\rfloor.$$
This means that
$$\sum_{m=0}^{j-1} \prod_{i=1}^a |\Lambda_{\ell_i, m}^+| = \sum_{m=0}^{j-1} \prod_{i=1}^a \Theta(m^{\left\lfloor \frac{\ell_i}2 \right\rfloor}) \leq \sum_{m=0}^{j-1} \Theta(m^{\left\lfloor \frac{n-2}2 \right\rfloor}) = \Theta(j^{\left\lfloor \frac{n}2 \right\rfloor}).$$
Therefore, $|\Lambda_{n,j}^+| \leq \Theta(j^{\left\lfloor \frac n2 \right\rfloor})$. We claim that equality is achieved. Consider $\alpha = (2, 1, \dots, 1)$. Then $\ell_1 = n-2$ and $\ell_2 = 1$, so $\left\lfloor \frac{\ell_1}2 \right\rfloor + \left\lfloor \frac{\ell_2}2 \right\rfloor =  \left\lfloor \frac{n-2}2 \right\rfloor + \left\lfloor \frac12 \right\rfloor = \left\lfloor \frac{n-2}2 \right\rfloor$, achieving equality. Therefore, $|\Lambda_{n,j}^+| = \Theta(j^{\left\lfloor \frac n2 \right\rfloor})$, concluding the induction.
\end{proof}

We can further refine this asymptotic by computing the leading coefficient in the following theorem.

\begin{theorem}\label{thm:asymptotic}
We have
$$|\Lambda_{n,k}^+| \sim \frac{a_{\left\lfloor \frac{n+1}2 \right\rfloor}}{\left\lfloor \frac n2 \right\rfloor !} k^{\left\lfloor \frac n2 \right\rfloor},$$
where $a_i$ is a sequence defined by $a_0 = a_1 = 1$ and $a_i = a_{i-1} + (i-1)a_{i-2}$.
\end{theorem}
\begin{remark}
The numbers $a_i$ are sometimes called the telephone numbers (OEIS~\cite[A000085]{oeis}). Some examples of its appearances include:
\begin{enumerate}
    \item The number of self-inverse permutations (involutions) on $n$ letters,
    \item The number of standard Young tableaux with $n$ cells,
    \item The sum of the degrees of the irreducible representations of the symmetric group $S_n$.
\end{enumerate}
\end{remark}
\begin{proof}
From Proposition~\ref{prop:lambda_nk_asymptotic}, we know that $|\Lambda_{n,k}^+| = \Theta(k^{\left\lfloor \frac n2 \right\rfloor})$. It remains to find the coefficient of the leading term in the recursion given by Theorem~\ref{thm:lambda_nk_recursion}. Let $b_n$ be the leading coefficient, i.e.~$|\Lambda_{n,k}^+| \sim b_n k^{\left\lfloor \frac n2 \right\rfloor}$. We know that $b_0, \dots, b_4$ are $1, 1, 1, 2, 1$, respectively.

As in the proof of Proposition~\ref{prop:lambda_nk_asymptotic}, the only terms that contribute to this leading term correspond to partitions $\alpha \vdash n$ such that $\alpha = (a^{\ell_a},\dots,1^{\ell_1})$ such that
$$\left\lfloor \frac{\ell_1}2 \right\rfloor + \cdots + \left\lfloor \frac{\ell_a}2 \right\rfloor = \left\lfloor \frac{n-2}2 \right\rfloor.$$
If $n$ is even, the only partitions that satisfy this are $(2, 1, \dots, 1)$ and $(2, 2, 1, \dots, 1)$. This means that
$$|\Lambda_{n,k}^+| \sim \sum_{m=0}^{k-1} b_1 \cdot b_{n-2} k^{\left\lfloor\frac n2 \right\rfloor -1} + \sum_{m=0}^{k-1} b_2 \cdot b_{n-4} k^{\left\lfloor\frac n2 \right\rfloor -1} \sim \frac{b_{n-2} + b_{n-4}}{\left\lfloor \frac n2 \right\rfloor} k^{\left\lfloor \frac n2 \right\rfloor},$$
where we use Faulhaber's formula for the sum of powers. Therefore, if $n$ is even, $b_n = \frac{2(b_{n-2} + b_{n-4})}{n}$.

On the other hand, if $n$ is odd, the only partitions that satisfy this are $(2, 1, \dots, 1), (3, 1, \dots, 1),$ and $(2, 2, 1, \dots, 1)$. This means that
$$|\Lambda_{n,k}^+| \sim \sum_{m=0}^{k-1} (b_1b_{n-2} + b_1b_{n-3} + b_2b_{n-4})k^{\left\lfloor \frac n2 \right\rfloor -1} \sim \frac{b_{n-2} + b_{n-3} + b_{n-4}}{\left\lfloor \frac n2 \right\rfloor} k^{\left\lfloor \frac n2 \right\rfloor}.$$
Therefore, if $n$ is odd, then $b_n = \frac{2(b_{n-2} + b_{n-3} + b_{n-4})}{n-1}$.

Now, we claim that $b_n = \frac{a_{\left\lfloor \frac{n+1}2 \right\rfloor}}{\left\lfloor \frac n2 \right\rfloor !}$ where $a_i$ is the sequence defined in the statement. This can be checked to work for the base cases. Using induction and by direct substitution, we could show that this also satisfies the two recursive equations we found above for $b_n$ for both even and odd $n$.
\end{proof}

\subsection{Understanding the asymptotics}\label{subsec:understanding-asymptotic}

The appearance of the telephone numbers is quite interesting, especially with its natural connections to permutations and partitions. In an effort to explain this relationship, we make the following note. First, each iteration brings the magnitude of the numbers down by a factor of $p$, so if we end up at zeroes after at most $k$ iterations, then the initial numbers must be at most on the order of $p^k$. Also, this means that each number is going to be on the order of a power of $p$, plus some smaller order terms. Furthermore, in Section~\ref{subsec:antisymmetry}, we show that each distinguished weight is anti-symmetric, so each sequence is identified by the first $\left\lfloor \frac n2 \right\rfloor$ numbers in the sequence. Putting these facts together, we find that each distinguished weight can be approximated by an anti-symmetric, weakly decreasing sequence of powers of $p$. Also, this implies our asymptotic also approximates the number of distinguished weights in the cube $[-p^k, p^k]^{n}$.

Furthermore, this intuition aligns with the asymptotics found in the above theorem; we can count the number of weakly decreasing sequences of length $\left\lfloor \frac n2 \right\rfloor$ such that each can be drawn from the set $\{p^0, \dots, p^{k-1}\}$. By stars and bars, this is $\binom{\left\lfloor \frac n2 \right\rfloor + k -1 }{k-1} = \frac{(\left\lfloor \frac n2 \right\rfloor + k -1) \cdots (k)}{\left\lfloor \frac n2 \right\rfloor !}$. When $k \gg n$, this is asymptotically $\frac{k^{\left\lfloor \frac n2 \right\rfloor}}{\left\lfloor \frac n2 \right\rfloor !}$, which is exactly the form in the theorem above without the telephone number coefficient.

For $n=4$, all distinguished weights are in the form $(x,y,-y,-x)$. Graphing the points $(x,y)$ for $p=5$ and at most $20$ iterations, we get the scatter plot shown below in Figure~\ref{fig:n4p5}. Similarly, for $n=5$, $p=11$, and at most $15$ iterations, all distinguished weights are in the form $(x,y,0,-y,-x)$, and the graph that arises is shown below in Figure~\ref{fig:n5p11}. First, we can observe a very regular pattern outside of the diagonals. To be precise, we can observe that for $(k_1, k_2)$ where $k_1 > k_2 + 2$, then there is $1$ distinguished weight $(x,y,-y,-x)$ and $3$ distinguished weights $(x,y,0,-y,-x)$ such that $p^{k_1} \leq x \leq p^{k_1+1}$ and $p^{k_2} \leq y \leq p^{k_2+1}$. Second, we can observe the quadratic growth that we proved in Proposition~\ref{prop:lambda_nk_asymptotic} through the ``interior" points. It should be noted that if we zoom in, there is also quadratic growth of the points along the diagonal $y=x$, which appear in very small and tight clusters and are almost unnoticeable in Figure~\ref{fig:n4-n5-examples}.

\begin{figure}[h]
    \centering
    \begin{subfigure}[H]{0.45\textwidth}
        \centering
        \includegraphics[width=2.8in]{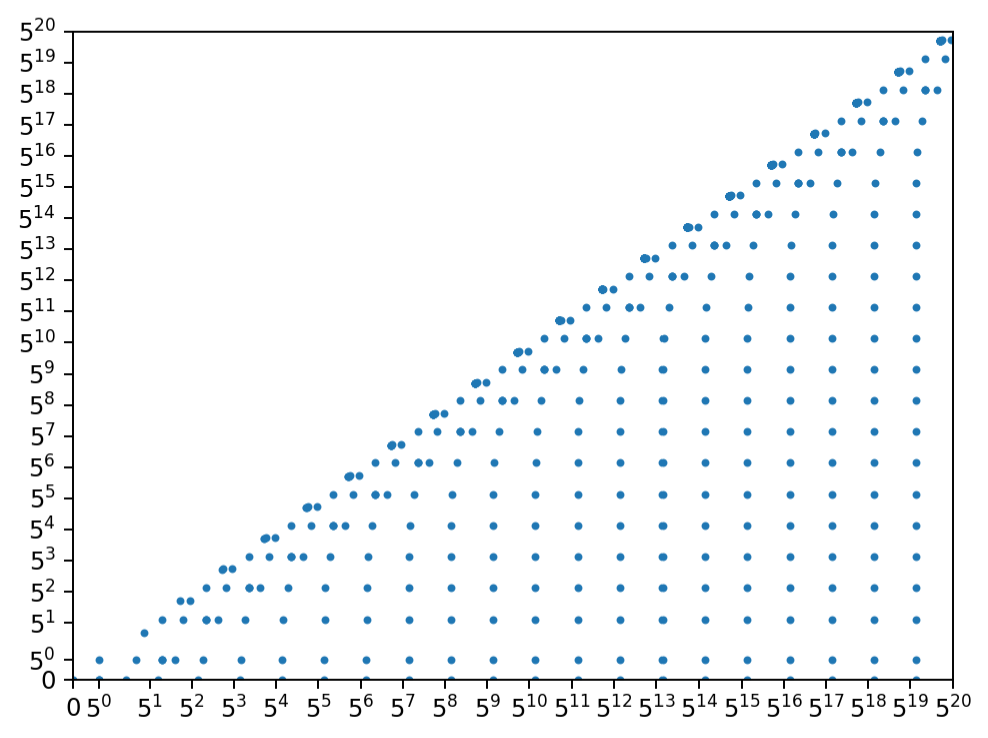}
        \caption{Scatter plot of $(x,y)$ for all distinguished weights $(x,y,-y,-x) \in \Lambda_4^+$ where $p=5$. Axes are log-scaled.}
        \label{fig:n4p5}
    \end{subfigure}
    \hfill
    \begin{subfigure}[H]{0.45\textwidth}
        \centering
        \includegraphics[width=2.83in]{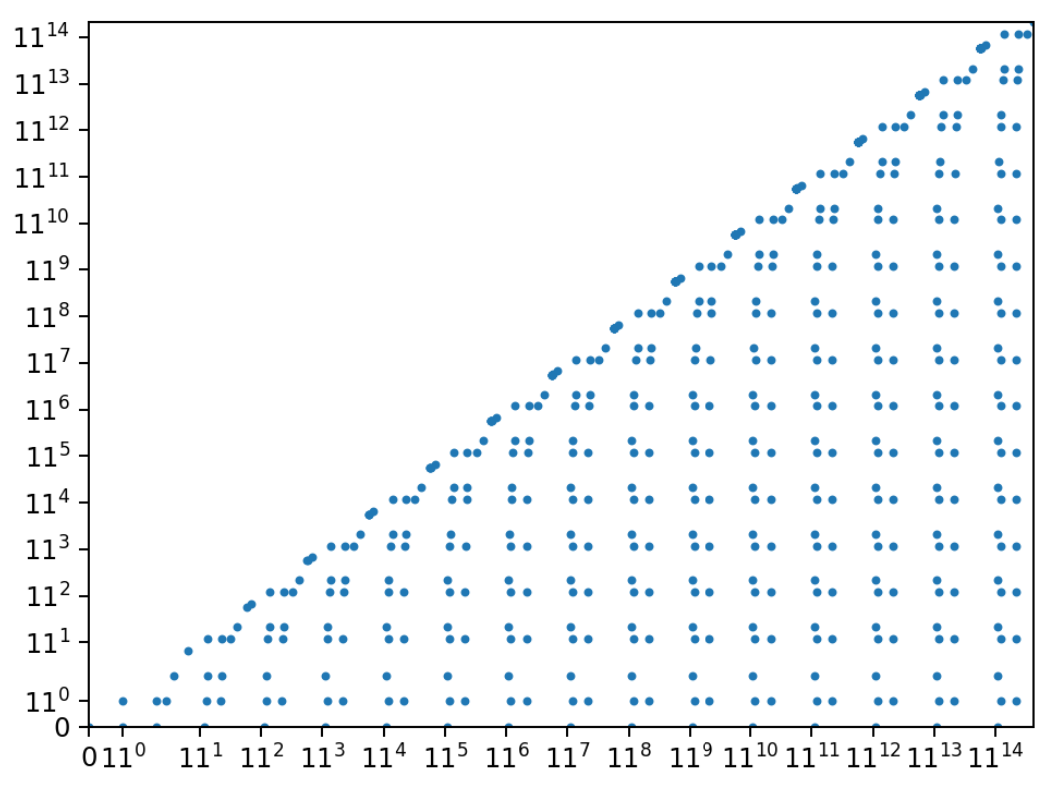}
        \caption{Scatter plot of $(x,y)$ for all distinguished weights $(x,y,0,-y,-x) \in \Lambda_5^+$ where $p=11$. Axes are log-scaled.}
        \label{fig:n5p11}
    \end{subfigure}
    \caption{}
    \label{fig:n4-n5-examples}
\end{figure}

This leads us to the following observation: suppose we have nonnegative integers $k_1, \dots, k_{\left\lfloor \frac n2 \right\rfloor}$ such that $k_i \geq k_{i+1} + \epsilon$ for $1 \leq i < \left\lfloor \frac n2 \right\rfloor$ for some $\epsilon > 0$ large enough. Then the number of distinguished weights $\lambda = (\lambda_1, \lambda_2, \dots, -\lambda_1)$ such that $p^{k_i} \leq \lambda_i < p^{k_i+1}$ for $1 \leq i \leq \left\lfloor \frac n2 \right\rfloor$ is $1$ if $n$ is even (as in the $n=4$ case above) and $\frac{n+1}2$ if $n$ is odd (as in the $n=5$ case above). Visually, this is stating that if $n$ is even, then there is one distinguished weight near each point with coordinates as powers of $p$, and if $n$ is odd, then there is $\frac{n+1}2$ distinguished weights near each point with coordinates as powers of $p$, assuming that we are sufficiently far from the diagonals, i.e. where two coordinates have similar powers of $p$. This can be unrigorously seen recursively:

If $n$ is even, then, as in the proof of Theorem~\ref{thm:asymptotic}, the leading term of the asymptotic comes recursively from the partitions $(2, 1, \dots, 1)$ and $(2, 2, 1, \dots, 1)$. We can explicitly describe the inverse of
$$((\lambda_1, \lambda_2, \dots, \lambda_{\frac n2 -1}, -\lambda_{\frac n2 -1}, \dots, -\lambda_1), (0))$$
and
$$((\lambda_1, \lambda_2, \dots, \lambda_{\frac n2 -2}, -\lambda_{\frac n2 -2}, \dots, -\lambda_1), (0,0))$$
through iterations of $LV_p^{-1}$ and keep track of just the largest order of power of $p$ in each term of the sequence. Then we could see inductively that if the $n-2$ case only has $1$ distinguished weight with a particular decreasing sequence of powers of $p$ with large enough consecutive differences, then corresponds directly in the case of $n$. Since there is only $1$ for $n=2$, then the count always is the same for even $n$.

If $n$ is odd, then similar analysis shows that there is always $1$ more than the $n-2$ case, so the number of these points grow linearly with $n$ as $\frac{n+1}2$.

However, this discussion does not explain combinatorially the appearance of the telephone numbers in the asymptotics of $|\Lambda_{n,k}^+|$. In particular, it explains the regularity in the ``interior" of the space of distinguished weights, but the points near the boundaries, i.e. where two coordinates are on or near the same order of power of $p$, are not explained. We would like to have a more direct way of interpreting these distinguished weights and the leading term of this asymptotic than through a recursive equation.

\newpage
\section*{Acknowledgements}


This research was conducted as a part of the MIT Summer Program for Undergraduate Research (SPUR). I am grateful to Andrei Ionov and Merrick Cai for their unwavering support and guidance throughout this research. I would like to thank Roman Bezrukavnikov for suggesting this project and to thank Vasily Krylov for his valuable communications. Additionally, I appreciate David Jerison and Jonathan Bloom for their thoughtful feedback throughout the project.

\bibliographystyle{plain}

\end{document}